\documentclass[draft]{article}
\def\today{13.09.10} 
\usepackage{amsmath,amsfonts,amsthm,amssymb,amscd}

\renewcommand{\Re}{\mathop{\rm Re}\nolimits}
\renewcommand{\Im}{\mathop{\rm Im}\nolimits}
\def\S{\mathhexbox278}

\newcommand{\beq}{\begin{equation}}
\newcommand{\ee}{\end{equation}}

\theoremstyle{plain} \newtheorem{theorem}{Theorem}[section]
\newtheorem{lemma}[theorem]{Lemma}
\newtheorem{proposition}[theorem]{Proposition}
 \theoremstyle{definition}
\newtheorem{definition}[theorem]{Definition} \theoremstyle{remark}
\newtheorem{remark}[theorem]{Remark}

\newcommand{\R}{{\mathbb R}} \newcommand{\U}{{\mathcal U}}

\newcommand{\Z}{{\mathbb Z}}
\newcommand{\Nn}{{\mathbb N}}

\newcommand{\E}{{\mathcal E}}

\newcommand{\Tr}{{\mathcal T}}
\newcommand{\Tra}{{\phi}}

\newcommand{\Ph}{{\mathcal P}}

\newcommand{\M}{{\mathcal M}}

\newcommand{\resto}{{\mathcal R}}
\def\im{{\rm i}}
\newcommand{\N}{{\mathcal N}}
\newcommand{\Nc}{{\mathcal N}}
\newcommand{\Sc}{{\mathcal S}}

\newcommand{\C}{\mathbb{C}}

\newcommand{\G}{{\mathcal G}}

\font\strana=cmti10
\def\lie{\hbox{\strana \char'44}}
\def\uno{{\kern+.3em {\rm 1} \kern -.22em {\rm l}}}
\def\beffe{{\bf f}}
\def\Phib{{\bf \Phi}}

\def\norma#1{\left\| #1\right\|}

\numberwithin{equation}{section}

\begin{document}

\title{On dispersion of
small energy solutions of the nonlinear Klein Gordon equation with a
potential}

\author {Dario Bambusi, Scipio Cuccagna}

\date{\today}
\maketitle

\begin{abstract}
In this paper we study small amplitude solutions of nonlinear Klein
Gordon equations with a potential. Under suitable smoothness and
decay assumptions on the potential  and a genericity assumption on
the nonlinearity, we prove that all small  energy solutions are
asymptotically free. In   cases where the linear system has at most
one bound state the result was already proved by Soffer and
Weinstein: we obtain here a result valid in the case of an arbitrary
number of possibly degenerate bound states. The proof is based on a
combination of Birkhoff normal form techniques and dispersive
estimates.
\end{abstract}

\section{Introduction}
\label{1} In this paper we study small amplitude solutions of the
nonlinear Klein Gordon equation (NLKG)
\begin{equation}\label{NLKG}
 u_{tt}-\Delta u +Vu+m^2u+\beta '(u) =0, \qquad (t,x)\in\mathbb{
   R}\times \mathbb{ R}^3
\end{equation}
with $ -\Delta +V(x)+m^2$ a positive short range Schr\"odinger
operator, and $\beta'$ a smooth function having a zero of order 3 at
the origin and growing at most like $u^3$ at infinity. Under
suitable smoothness and decay properties on the potential $V$ and on
$\beta '$, and under a genericity assumption on the nonlinearity,  to
be discussed below, we prove that all small  energy solutions
are asymptotically free. Thus in particular the system does not
admit  small energy periodic or quasiperiodic solutions,   in
contrast  with what happens in bounded domains where KAM theory can
be used to prove existence of quasiperiodic solutions
\cite{K1,CW92,W90,Bou03,EK06}.

A crucial role in our discussion   is played by the spectrum of the
Schr\"odinger operator $ -\Delta +V(x)$. If $ -\Delta +V(x)$
 does not have eigenvalues, then the asymptotic freedom of
solutions follows from a perturbative argument based on a theorem by
Yajima \cite{yajima}. If
   $-\Delta +V
+m^2$ has just one nondegenerate eigenvalue lying close to the
continuous spectrum, then the result is proved by
\cite{sofferweinsten1}. We generalize this result,   easing most
restrictions on the spectrum of $-\Delta +V +m^2$.

From a technical standpoint, the key   is to prove that, due to
nonlinear coupling, there is leaking of energy from the discrete
modes to the continuous ones. The continuous modes should disperse
  by perturbation,   because of the linear dispersion.  In
\cite{sofferweinsten1} this leaking occurs because the discrete mode
equation has a key coefficient of  positive   sign, which yields
dissipation. In \cite{sofferweinsten1} this coefficient is of the
form $\langle D F,F\rangle $ for $D $ a
  positive operator and $F$ a  function. Assuming the
generic condition $\langle D F,F\rangle\not=0 $ (which is called
nonlinear Fermi golden rule or FGR), then such a quantity is
strictly positive. This gives rise to dissipative effects leading to
the result. The presence of terms of the form $\langle D F,F\rangle
$ was first pointed out and exploited for nonlinear problems in
\cite{sigal}, which proves that periodic and quasiperiodic solutions
of the linear equation are unstable with respect to nonlinear
perturbations. In the problem treated in \cite{sigal}, this
coefficient appears directly. In our case, to   exploit the
coefficient   it is first necessary to simplify the equations by
means of  normal form expansions. The normal forms argument was
first introduced in \cite{busper2}, later by \cite{sofferweinsten1},
(see also \cite{zhousigal,cuccagnamizumachi} and for further
references \cite{cuccagnatarulli}).

In the case    when the  eigenvalues of $-\Delta +V +m^2$ are not
close to the continuous spectrum,   the crucial coefficients in the
equations of the discrete modes are of the form $\langle D
F,G\rangle $ for $F$ and $G$ not obviously related, if one follows
the scheme in
\cite{busper2,sofferweinsten1,zhousigal,cuccagnamizumachi}. The
argument in \cite{cuccagnamizumachi} shows indirectly that, in the
case of just one simple eigenvalue, this coefficient is semidefinite
positive. But this is not clear any more in the case of multiple
eigenvalues of possibly high multiplicity, if one follows the scheme
in \cite{busper2,sofferweinsten1,zhousigal,cuccagnamizumachi}. In
the present paper we fill this gap. Using the Hamiltonian structure
of (\ref{NLKG})  and the Birkhoff normal form theory, we show that
dissipativity is a generic feature of the problem. Here lies the
novelty of this paper:  previous references  perform normal form
expansions losing sight of the Hamiltonian structure of
(\ref{NLKG}). It turns out that the Hamiltonian   structure is
crucial.

We recall that Birkhoff normal form theory has been recently
extended to a quite large class of Hamiltonian partial differential
equations (see for example \cite{BN98,Bam03,BG06}). However here we
 need  to deal with two specific issues. The
first one is that we need  to produce a normal form which keeps some
memory of the fact that the original Hamiltonian is local, since
locality is a fundamental property needed for the dispersive
estimates used   to prove dissipation. The second issue is that the
Hamiltonian function (and its vector field) of the NLKG has only
finite regularity, so it is not a priori obvious how to put the
system in normal form at high order. This problem is here solved by
noticing  that our normal form  is needed  only to simplify the
dependence on the discrete modes and to decouple the discrete modes
from the continuous ones. This can be obtained by a coherent
recursive construction yielding   analytic canonical
transformations.

Finally, the related problem of asymptotic stability of ground states
of the NLS initiated in \cite{sofferweinsten2}, see also the seminal
papers \cite{sofferweinsten3,busper1,busper2,zhousigal}, has been
solved in \cite{cu1} drawing the ideas in the present paper. Other
references on the NLS which we mention later are
\cite{tsai,zhouweinstein1}.  For further references we refer to
\cite{cuccagnatarulli,cu1}.

\section{Statement of the main result}
\label{statement}

We begin by stating our assumptions.

\begin{itemize}
\item[(H1)] $V(x)$ is real valued and $|\partial ^\alpha _x V(x)|\le C
  \langle x\rangle ^{-5-\sigma}$ for $|\alpha |\le 2$, where $C>0$ and
  $\sigma >0$ are fixed constants and $\langle x\rangle :=
  \sqrt{1+|x|^2} $; $V(x)$ is smooth with $|\partial ^\alpha _x
  V(x)|\le C_\alpha <\infty$ for all $\alpha$;

\item[(H2)] 0 is  neither  an eigenvalue nor a resonance for $ -\Delta +V
$, i.e. there are no   nonzero solutions of $\Delta u=Vu$ in $
\mathbb{R}^3$ with $|u(x)|\lesssim \langle x\rangle ^{-1}.$
\end{itemize}
It is well known that (H1)--(H2) imply that the set of eigenvalues
$\sigma _d(-\Delta +V)\equiv\left\{-\lambda_j^2\right\}_{j=1}^n$ is
finite, contained in $(-\infty ,0)$, with each eigenvalue of finite
multiplicity. We take a mass term $m^2$ such that $-\Delta +V +m^2> 0$
and we assume that indexes have been chosen so that $-\lambda _1^2\le
\dots \le -\lambda _n^2$. We set $\omega _j=\omega _j(m):=\sqrt{ m^2-
  \lambda _j^2}$. We assume $m>0$ and $\lambda _j>0$.
Notice that the $ \lambda _j $ are not necessarily pairwise
distinct. We assume that $m $ is not a multiple of any of the
$\omega _j$'s:

\begin{itemize}
\item[(H3)]   for any $\omega _j$   there exists    $N_j\in \mathbb{N}$
 such that $N_j\omega _j<m<  (N_j+1)\omega _j$.
\end{itemize}
Notice that  $N_1=N:=\sup _jN_j$. Hypothesis (H3) is a special case
of the following  hypothesis:

\begin{itemize}
\item[(H4)]   there is no multi index $\mu \in \mathbb{Z}^{n}$
with $|\mu|:=|\mu_1|+...+|\mu_n|\leq 2N_1+3$ such that $\mu \cdot
\omega =m$.
\end{itemize}
 We furthermore require:

\begin{itemize}
\item[(H5)] if $\omega _{j_1}<...<\omega _{j_k}$ are $k$ distinct
  $\omega$'s, and $\mu\in \Z^k$ satisfies
  $|\mu| \leq 2N_1+3$, then we have
$$
\mu _1\omega _{j_1}+\dots +\mu _k\omega
_{j_k}=0 \iff \mu=0\ .
$$
\end{itemize}
\begin{remark}
\label{rem.in.1}  Using the fact that for any
  $\mu$ the quantities $\mu \cdot \omega $ are holomorphic functions
    in $m$ for $\Re m > \lambda_1$, it is easy to show that there exists a
discrete set $D\subset ( \lambda_1 ,\infty)$, such that for $m\not \in
D$ hypotheses (H3-H5) are true.
\end{remark}

Assumptions (H1)--(H5) refer to the properties of the linear part of
the equation. Consider now $\beta (u)=\int _0^u\beta '(s) ds.$ We
assume the following hypothesis:
\begin{itemize}
\item[(H6)]    we assume
that there
  exists a smooth function $\tilde \beta\in C^{\infty}(\R,\R)$ such
  that $\beta(u)=u^4\tilde \beta(u)$ and, for any $j\geq 0$ there exists
$C_j>0$ such that $|\tilde \beta ^{(j)}(u)|\le C_j\langle u\rangle^{-j}$.
\end{itemize}

Finally there is an hypothesis relating the linear operator $-\Delta
+V+m^2$ and the nonlinearity $\beta (u)$. It is a nondegeneracy
hypothesis that, following \cite{sigal,sofferweinsten1}, we call
nonlinear Fermi golden rule. Specifically, the main result of this
paper is that certain coefficients  related to
the resonance between discrete and continuous modes are non
negative. The nondegeneracy hypothesis is that they are strictly
positive. We show in Proposition \ref{H6} that this hypothesis holds
generically, in some sense. The precise statement of the  hypothesis
requires some notation  and preliminaries, so   is deferred to
section \ref{subsection:nondegeneracy}. We assume what follows:
\begin{itemize}
\item[(H7)] we assume that \eqref{m.34} or, equivalently
  \eqref{m.34a}, holds.
\end{itemize}

  (H7) is the most significant of our hypotheses. It should hold
  quite generally. By way of illustration,
in Section \ref{subsection:nondegeneracy} we   prove the following
result:
\begin{proposition}
\label{H6} Assume that $V$ satisfies (H1)--(H2), decreases
exponentially  together  with all its derivatives as
$|x|\to\infty$ and that all the eigenvalues of $-\Delta +V$ are
simple.  Then there exist a finite set $\M\subset( \lambda_1 ,+\infty
)$, for any $m\in ( \lambda_1 ,+\infty )\backslash\M$ a finite set
$\widehat{M} (m)\subset\Z ^n$ locally constant in $m$, functions
$f_{\mu ,m}^{(\pm )} \in C^\infty ( \R ^{|\mu |-4}, \R )$ for $\mu\in
\widehat{M} (m)$, such that (H7) holds if the following is true: $m\in
( \lambda_1 ,+\infty )\backslash\M$ and for both signs $\pm$
$$
\beta _{ |\mu | } \not=f_{\mu ,m}^{(\pm )} ( \beta _4,...,\beta _{
|\mu |-1 } ) \text{ for all $\mu \in \widehat{M} (m)$ and where
$\beta _j:=\beta ^{(j)}(0)/j!$.}
$$
\end{proposition}
Now we state the main result of this paper. Denote
$K_0(t)=\frac{\sin (t \sqrt{-\Delta +m^2})}{\sqrt{-\Delta +m^2}}.$
Then we prove:
\begin{theorem}\label{theorem-1.1}
Assume hypotheses (H1)--(H7). Then there exist $\varepsilon _0>0$
and $C>0$ such that for any $\| (u_0,v_0) \| _{H^1\times L^2}\le
\epsilon <\varepsilon _0$  the solution of (\ref{NLKG}) with
$(u(0),u_t(0))=(u_0,v_0)$ is globally defined and there are $(u_\pm
,v_\pm )$ with $\| (u_\pm ,v_\pm ) \| _{H^1\times L^2}\le C \epsilon
$
\begin{equation}\label{scattering}\lim _{t\to \pm \infty } \|  u(t) -
 K'_0(t)
u_\pm - K_0(t) v_\pm \| _{H^1 }=0.
\end{equation}
It is possible to write $u(t,x)=A(t,x)+\widetilde{u}(t,x)$  with
$|A(t,x)|\le C_N(t) \langle x \rangle ^{-N}$ for any $N$, with $\lim
_{|t|\to \infty }C_N(t)=0$ and such that for any pair $(r,p)$ which
is admissible, by which we mean that
\begin{equation}\label{admissiblepair}  2/r+3/p= 3/2\,
 , \quad 6\ge p\ge 2\, , \quad
r\ge 2,
\end{equation}
we have
\begin{equation}\label{Strichartz} \|  \widetilde{u} \| _{L^r_tW^{
\frac{1}{p}-\frac{1}{r}+\frac{1}{2},p}_x}\le
 C\| (u_0,v_0) \| _{H^1\times L^2}.
\end{equation}
\end{theorem}
\begin{remark}Theorem \ref{theorem-1.1} is well known in the particular case
$V=0$, see Theorem 6.2.1 \cite{cazenave}.  In this case $
\widetilde{u}=u$. If the operator $-\Delta +V $
  does not have eigenvalues and satisfies the estimates in Lemma
  \ref{lemma-strichartzflat}, then Theorem \ref{theorem-1.1} continues
  to hold. Work by Yajima \cite{yajima} guarantees that this indeed is
  the case for operators satisfying (H1)--(H2) such that $\sigma
  _d(-\Delta +V)$ is empty, see Lemma \ref{lemma-strichartz1}. These
  results are obtained  by  thinking the nonlinear problem as a
  perturbation of the linear problem.
\end{remark}

\begin{remark} Theorem \ref{theorem-1.1} can be thought as an
asymptotic  stability result of the 0 solution. Stability is well
known, see Theorem \ref{theorem-globalwellposedness} below.
\end{remark}

\begin{remark} Theorem \ref{theorem-1.1}  in the case  when $\sigma
  _d(-\Delta +V)$ consists of a single eigenvalue can be proved
  following a simpler version of the argument in
  \cite{cuccagnamizumachi}.
\end{remark}
\begin{remark}
Theorem \ref{theorem-1.1}  in the case  when $\sigma _d(-\Delta +V)$
consists of a single eigenvalue $-\lambda ^2$ such that for $\omega
=\sqrt{m^2-\lambda ^2}$ we have $3\omega >m$ is proved in
\cite{sofferweinsten1} assuming $\| (u_0,v_0) \| _{(H^2\cap
W^{2,1})\times (H^1\cap W^{1,1})}$ small. Notice that formula (1.10)
\cite{sofferweinsten1} contains a decay rate of dispersion of the
various components of $ {u}(t)$. For the initial data in the larger
class considered in Theorem \ref{theorem-1.1}, such kind of decay
rates cannot be proved. Restricting  initial data to the class in
\cite{sofferweinsten1}, it is possible to prove appropriate decay
rates also for the solutions in Theorem \ref{theorem-1.1}.
\end{remark}

\begin{remark}   Theorem \ref{theorem-1.1}
is stated only for $\R ^d$ with $d=3$.  Versions of this theorem can
be proved for any $d$. In particular, the crux of the paper, that is
the normal form expansion in Theorem \ref{main} and the discussion
of the discrete modes, are not affected by the spatial dimension.
\end{remark}
In view of the above remarks,
   we focus our
attention to the case when $-\Delta +V$ admits eigenvalues,
especially the case of many eigenvalues.

We end this section with some notation. Given two functions
$f,g:\mathbb{R}^3\to \mathbb{C}$ we set $\langle f,g\rangle = \int
_{\mathbb{R}^3}f(x) g(x) dx$. For $k\in \mathbb{R}$ and $1< p< \infty
$ we denote for $K=\mathbb{R},\mathbb{C}$
 \[ W^{k,p}(\mathbb{R}^3,K)=\{ f:\mathbb{R}^3\to K \text{ s.t.}
 \|   f\| _{W^{k,p}}:=\| (-\Delta +1)^{k/2} f\| _{L^p} <\infty .\}\]
In particular we set $H^k(\mathbb{R}^3,K)=W^{k,2}(\mathbb{R}^3,K) $
and $L^p(\mathbb{R}^3,K)=W^{0,p}(\mathbb{R}^3,K) $. For $p=1,\infty$
and $k\in \mathbb{N}$ we denote by $W^{k,p}(\mathbb{R}^3,K)$ the
functions such that  $\partial _x^\alpha f\in
 L^{ p}(\mathbb{R}^3,K)$ for all $|\alpha |\le k$  (we recall that for
 $1<p<\infty$ the two definitions of  $W^{k,p}$ yield the same space).
For any $s\in \mathbb{R}$ we set

 \[ H^{ k,s}(\mathbb{R}^3,K)=\{ f:\mathbb{R}^3\to K \text{ s.t.}
 \| f\| _{H^{s,k}}:=\| \langle x \rangle ^s(-\Delta +1)^{k/2} f\| _{L^2
 }<\infty \}.\] In particular we set
 $L^{2,s}(\mathbb{R}^3,K)=H^{0,s}(\mathbb{R}^3,K).$
Sometimes, to emphasize that these spaces refer to spatial
variables, we will denote them by $W^{k,p}_x$, $L^{ p}_x$, $H^k_x$,
$H^{ k,s}_x$ and $L^{2,s}_x$. For $I$ an interval and $Y_x$ any of
these spaces, we will consider Banach spaces $L^p_t( I, Y_x)$ with
mixed norm $ \| f\| _{L^p_t( I, Y_x)}:= \| \| f\| _{Y_x} \| _{L^p_t(
I )}.$ Given an operator $A$, we will denote by $R_A(z)=(A-z)^{-1}$
its resolvent. We set $\mathbb{N}_0=\mathbb{N}\cup \{0 \}$. We will
consider multi indexes $\mu \in \mathbb{N}_0^n$. For   $\mu \in
\mathbb{Z}^n$ with $\mu =(\mu _1,..., \mu _n)$ we set $|\mu |=\sum
_{j=1}^n |\mu _j|.$ We also consider the set of Schwartz functions
$\mathcal{S}(\R ^3,\C )$ whose elements are the functions  $f\in
C^\infty (\R ^3,\C ) $ such that $\langle x \rangle ^N \partial
_x^\alpha f(x)\in L^\infty (\R ^3  )$ for all $N\in   \mathbb{{N}}
\cup \{0\} $ and $\alpha\in ( \mathbb{N} \cup \{0\})^3$.

\section{Global well posedness and Hamiltonian structure}
\label{2}

In $H^1(\R^3,\R)\times L^2(\R^3,\R)$ endowed with the standard
symplectic form, namely
\begin{eqnarray}
\label{symplecticform} \Omega((u_1,v_1);(u_2,v_2) ):=\langle
u_1,v_2\rangle_{L^2}- \langle u_2,v_1\rangle_{L^2}
\end{eqnarray}
we consider the Hamiltonian
\begin{eqnarray}
\label{H.1} H=H_L+H_P\ ,
\\
H_L:=\int_{\R^3}\frac{1}{2}(v^2+|\nabla u|^2+Vu^2+m^2u^2)d x\ ,
\nonumber \\
H_P:= \int_{\R^3}\beta (u )d x. \nonumber
\end{eqnarray}
The corresponding Hamilton equations are $\dot v=-\nabla_u H$, $\dot
u=\nabla_v H$, where $\nabla _uH$ is the gradient with respect to
the $L^2$ metric, explicitly defined by
$$
\langle\nabla_u H(u),h\rangle= d_u H(u) h\ ,\quad \forall h\in H^1\ ,
$$ and $d_u H(u)$ is the Frech\'et derivative of $H$ with respect to
$u$. It is easy to see that the Hamilton equations are explicitly
given by
\begin{equation}
\label{Hamiltonequations} \left(\dot v= \Delta u-Vu-m^2u-\beta '(u)\
,\, \dot u= v\right)\iff \ddot u= \Delta u-Vu-m^2u-\beta '(u)
\end{equation}
  First we recall that the NLKG (\ref{NLKG}) is globally
well posed for small initial data.
\begin{theorem}\label{theorem-globalwellposedness}
Assume $V\in L^{ p}_x$ with $p>3/2$. Then there exist $\varepsilon
_0>0$ and $C>0$ such that for any $\| (u_0,v_0) \| _{H^1_x\times
  L^2_x}\le \epsilon <\varepsilon _0$ and if we set $v(t)=u_t(t)$ and
$v_0=u_t(0)$, equation (\ref{NLKG}) admits exactly one    solution

\begin{equation}\label{strongsolution}
u \in C^{0}(\mathbb{R}, H^1_x) \cap C^{1}(\mathbb{R},
L^2_x) \end{equation}
such that
$(u(0),v(0))=(u_0,v_0)$.
 The map $(u_0,v_0)\to (u(t),v(t))$ is
continuous from the ball     $\| (u_0,v_0) \| _{H^1_x\times L^2_x}
<\varepsilon _0$ to $C^{0}(I, H^1_x) \times C^{0}(I, L^2_x)$ for any
bounded interval $I$. The Hamiltonian $H(u(t),v(t)) $ is constant,
and
\begin{equation}\label{boundedenergynorm}
\| (u(t),v(t)) \| _{H^1_x\times L^2_x}\le C   \| (u_0,v_0) \|
_{H^1_x\times L^2_x}.\end{equation} We have the equality

\begin{equation}\label{duhamelflat}
 u(t) = K'_0(t)u_0+K _0(t)v_0-
 \int _0^t K _0(t-s) (Vu(s)+\beta '(u(s)))  ds.\end{equation}

\end{theorem}
For statement and  proof see \S 6.2 and 6.3 \cite{cazenaveharaux}.

We associate to any $-\lambda _j^2$ an $L^2$ eigenvector
$\varphi_j(x)$, real valued and normalized.  We have $\varphi_j\in
 \mathcal{S}(\mathbb{R}^3, \mathbb{R})$.   Set $P_d
u=\sum \langle u,\varphi _j\rangle \varphi _j$   and set
$P_c=1-P_d$, the projector    in $L^2$ associated to the continuous
spectrum. Denote
\begin{equation}
\label{e.1} u=\sum_jq_j\varphi_j+P_cu\ ,\quad
v=\sum_jp_j\varphi_j+P_cv .
\end{equation}
   We have
\begin{eqnarray}
\label{H.2} H_P=  \int_{\R^3} \beta \left (\sum_jq_j\varphi_j+P_cu
\right ) d x .
\end{eqnarray}
Introduce the operator
\begin{equation}
\label{e.2} B:=P_c(-\Delta +V +m^2)^{1/2}P_c \ ,
\end{equation}
and the complex variables
\begin{equation}
\label{e.8} \xi_j:=\frac{q_j\sqrt{\omega_j}+\im
  \frac{p_j}{\sqrt{\omega_j}}} {\sqrt2}\ ,\quad f:=\frac{B^{1/2}P_c u +\im
B^{-1/2} P_cv}{\sqrt2}\ .
\end{equation}
By Theorem \ref{yajima},   \eqref{e.8}
  defines an isomorphism between
$H^1(\R^3,\R)\times L^2(\R^3,\R)$ and $\Ph^{1/2,0}:=\C^n\oplus P_c
H^{1/2,0}(\R^3,\C)$, which {\it from now on will be our phase
space}. We will often represent functions (and maps) on the phase
space as functions of the variables $\xi_j,\bar\xi_j,f,\bar f$. By
this we mean that a function $F(\xi,\bar\xi,f,\bar f)$ is   the
composition of the maps
$$
(\xi,f)\mapsto (\xi,\bar\xi,f,\bar f)\mapsto
F(\xi,\bar\xi,f,\bar f)\ .
$$ Correspondingly we define $ \partial _{\xi_j} =\frac{1}{2} (
\partial _{\Re \xi _j}-\im\partial _{\Im \xi _j} ) $ and $ \partial
_{\bar\xi _j} =\frac{1}{2} ( \partial _{\Re \xi _j}+\im \partial _{\Im
\xi _j} ) $, and analogously $\nabla_f:= \frac{1}{2} ( \nabla _{\Re
f}-\im \nabla _{\Im f })$, $\nabla_{\bar f}:= \frac{1}{2} ( \nabla
_{\Re f}+\im\nabla _{\Im f} )$.

In terms of these variables the symplectic form has the form
\begin{eqnarray}
\label{e.bis}
\Omega((\xi^{(1)},f^{(1)});(\xi^{(2)},f^{(2)}))= 2\Re \left[ \im
  \left(\sum_{j}\xi_j^{(1)}\bar \xi_j^{(2)}+\langle f^{(1)},\bar
  f^{(2)}\rangle \right) \right] \\ \nonumber
=-\im \sum_{j}\left(\bar \xi_j^{(1)} \xi_j^{(2)}-\xi_j^{(1)}\bar
\xi_j^{(2)} \right) -\im \left ( \langle f^{(2)},\bar f^{(1)}\rangle - \langle
f^{(1)},\bar f^{(2)}\rangle \right )
\end{eqnarray}
and
the Hamilton equations take the form
\begin{equation}
\label{e.9} \dot\xi_j=-\im\frac{\partial H}{\partial\bar\xi_j}\
,\quad \dot f= -\im \nabla_{\bar f}H\ .
\end{equation}
The Hamiltonian vector field $X_H$ of a function is given by
\begin{equation}
\label{e.9a} X_H(\xi,\bar\xi,f,\bar f)=\left( -\im\frac{\partial
  H}{\partial\bar\xi}, \im\frac{\partial H}{\partial\xi},-\im
  \nabla_{\bar f}H, \im \nabla_{f}H  \right)
\end{equation}
We consider the Poisson bracket
\begin{equation}
\label{e.9b} \left\{H,K\right\}:= \im\sum_{j}\left(\frac{\partial
H}{\partial
  \xi_j}\frac{\partial K}{\partial\overline{\xi}_j}-
  \frac{\partial H}{\partial
  \overline{\xi}_j}\frac{\partial K}{\partial \xi_j} \right)+ \im \left\langle
  \nabla_{  f}H,\nabla _{\overline{f}}K \right\rangle-\im \left\langle
 \nabla_{\overline{f}} H,  \nabla_{  f}K  \right\rangle .
\end{equation}
We emphasize that if $H$ and $K$ are real valued, then $
\left\{H,K\right\}$ is real valued. Later we will consider
Hamiltonians for which \eqref{e.9b} makes sense.

We introduce now some further notations that we will use in the
 sequel.

\begin{itemize}

\item
We denote the phase spaces $\Ph^{k,s} =\C^n\times P_cH^{k,s}
(\mathbb{R}^3, \mathbb{C}) $ with the spectral decomposition
associated to $-\Delta +V$.

\item $\beffe:=(f,\bar f)$, and we will denote by $\Phib:=(\Phi,\Psi)$
  a pair of functions each of which is in $ \mathcal{S}(\R ^3, \C )$.
\item Given $\mu\in \Nn^{n}$ we denote $\xi^\mu:=
  \prod_{j}{\xi_j^{\mu_j}}$, and similarly for $\bar\xi^\nu$.

\item A point of the phase space will usually be denoted by
  $z\equiv(\xi,f)$.
\end{itemize}

The form of $H_L$ and of $H_P$ are respectively
\begin{eqnarray}
\label{e.10} H_L=\sum_{j=1}^{n}\omega_j|\xi_j|^2+\langle \bar f,Bf
\rangle .
\\
\label{e.10a}
H_P(\xi,f)=\int_{\R^3}\beta(\sum
 \frac{\xi _j+\bar \xi _j}{\sqrt{2\omega _j}}\varphi _j(x)+U(x) )dx
\end{eqnarray}
where we wrote for simplicity $U=B^{-\frac{1}{2}}(f+\bar
 f)/\sqrt{2}\equiv P_cu$.

We will need something more about  the nonlinearity. Consider the
Taylor expansion
\begin{equation}  {\beta
}( \sum \frac{\xi _j+\bar \xi _j}{\sqrt{2\omega _j}}
 \varphi _j+U)  =  \sum _{l=0}^{3}F_l (x, \xi  )  U ^{l} +
 F_4 (x,\xi, U) U ^4 \nonumber
\end{equation}
with
\begin{eqnarray} &  \label{taylor1}F_l (x, \xi   )= \frac{1}{l!}
 \beta ^{(l)} ( \sum \frac{\xi _j+\bar \xi _j}{\sqrt{2\omega _j}}
 \varphi _j ) \ ,\quad l=0,1,2,3 \\& \label{taylor2} F_4 (x, \xi  , U)=
 \int _0^1 \frac{(1-\tau )^{3}}{3!} \beta ^{(4)}( \sum \frac{\xi
 _j+\bar \xi _j}{\sqrt{2\omega _j}}\varphi _j +\tau U )d\tau .
\end{eqnarray}

\begin{lemma}\label{lemma:H_P} The following holds true. \begin{itemize}
\item[(1)] For $l\leq 3$, the    functions
  $\xi \to F_{l } (\cdot , \xi  )  $ are in
  $C^\infty ( \mathbb{C}^{n  },H^{k,s})$
  for any $k,s$,
  and
\begin{eqnarray}
  H_l(\xi , U)
     =
\int_{\R^3}
  F_{l} (x, \xi    ) U^l dx \nonumber
\end{eqnarray}
are $H_l\in C^\infty (\mathbb{C}^{n  }\times H^{1}, \R) $.  In
particular we have
  derivatives, for $\ell\leq l$,
 \begin{equation}
  \partial ^\alpha _\xi d^\ell _U H_l \left [ \otimes _{j=1}^{\ell}
 g_j\right ] =l  \cdots (l-\ell +1)
 \int_{\R^3}
   \partial ^\alpha _\xi F_{l} (x, \xi    ) U ^{l-\ell }(x) \prod _{j=1}^{\ell}
  g_j(x)  dx. \nonumber
 \end{equation}

 \item[(2)]   $F _{l } $ has a 0 of order
 $4-l$ at $\xi =0$:
$$
\norma{F_l(\cdot , \xi  )}_{H^{k,s}}\leq C\norma{\xi}^{4-l}.
$$

\item[(3)]  The map $\C^n\times\R^3\times \R \ni( \xi , x,    Y)\mapsto
  F_4( x, \xi   , Y)\in \R$ is $C^\infty$;   for any $k>0$ there
  exists $C_k$ such that
$ |\partial^k_Y F_4( x, \xi   ,Y)|\leq C_k\ . $  Denote
\begin{eqnarray}
  H_4(\xi,U) =
\int_{\R^3}
  F_{4} ( x, \xi    ,U(x))U^4(x) dx .\nonumber
\end{eqnarray}
Then the map $\C^n\ni\xi\mapsto H_4(\xi,.)\in C^2(H^1,\C )$ is
 $C^\infty$. In particular
\begin{eqnarray}
 \partial _\xi
  ^\alpha d_U H_4 [g ] =
\int_{\R^3} \partial _\xi
  ^\alpha \partial _Y
  \Psi  ( x, \xi   ,U(x) )   g (x) dx \nonumber
\end{eqnarray} where
$\Psi  ( x, \xi   ,Y)=F_{4} ( x, \xi   ,Y)Y ^4$.
\end{itemize}
\end{lemma}
\proof The result follows by standard computations and explicit
estimates of the remainder, see p. 59 \cite{cazenave}.\qed

\section{Normal form}

\subsection{Lie transform}

We will iteratively eliminate from the Hamiltonian monomials,
simplifying the part linear in $f$ and $\bar f$  and the part
  independent of such variables.  We will
use canonical transformations generated by Lie transform, namely
the time 1 flow of a suitable auxiliary Hamiltonian function.
Consider a function  $\chi$ of the form
\begin{equation}
\label{chi.1} \chi(z)\equiv
\chi(\xi,f)=\chi_0(\xi,\bar\xi)+\sum_{|\mu|+|\nu|=M_0 +1}
\xi^{\mu}\bar\xi^{\nu}\int_{\R^3}\Phib_{\mu,\nu}\cdot \beffe dx
\end{equation}
where $\Phib_{\mu,\nu}\cdot \beffe:= \Phi_{\mu,\nu} f+
{\Psi}_{\mu,\nu}\bar f$ with $\Phi _{\mu,\nu} , \Psi _{\mu,\nu} \in
\mathcal{S}(\mathbb{R}^3, \mathbb{C})$ and where $\chi_0$ is a
homogeneous polynomial of degree $M_0 +2$.  The Hamiltonian vector
field satisfies $X_\chi \in C^\infty (\Ph ^{-\kappa ,-s }, \Ph ^{ k ,
  \tau } )$ for any $k,\kappa,s,\tau\geq 0$. Moreover we have
\begin{equation}
\label{Lipschitz} \| X_{\chi}(z) \|  _{\Ph  ^{k,\tau }}
 \le C_{k,s,\kappa,\tau}\| z\| _{\Ph ^{-\kappa ,-s }}^{M_0 +1}
.\end{equation} Since $X_{\chi}$ is a smooth polynomial it is also
analytic.  Denote by $\phi^t$  the   flow generated by $X_{\chi}$. For
fixed $\kappa,s$, $\phi^t$ is defined in $\Ph ^{-\kappa ,-s }$ up to
any fixed time $\bar t$, in a sufficiently small neighborhood
$\mathcal{U}^{-\kappa,-s}$ of the origin. For $\Ph ^{ k , \tau
}\hookrightarrow \Ph ^{-\kappa ,-s }$, by \eqref{Lipschitz} the flow
$\phi^t$ is defined for $0\le t \le \bar t$ in
$\mathcal{U}^{-\kappa,-s} \cap \Ph ^{ k , \tau }$. Set
$\phi:=\phi^1\equiv \phi^t\big|_{t=1}$

\begin{definition}
\label{Lie} The canonical transformation $\phi$ will be called the
{\it Lie transform} generated by $\chi$.\qed
\end{definition}

\begin{remark}
\label{rem.re}
The function $\chi$ extends to an analytic function on the
complexification of the phase space, namely the space in which $\xi$
is independent of $\bar \xi$ and $f$ is independent of $\bar f$. If
the original function $\chi$ is real valued (as in our situation),
then $\chi$ takes real values when $f$ is the complex conjugated of
$\bar f$ and $\xi$ the complex conjugated of $\bar \xi$. In this case,
by the very construction, the Lie transform generated by $\chi$ leaves
invariant the submanifold of the complexified phase space
corresponding to the original real phase space.
\end{remark}

\begin{lemma}\label{lie_trans}
Consider a functional $\chi$ of the form \eqref{chi.1}. Assume
$\Phi_{\mu,\nu}, \Psi_{\mu,\nu}\in \mathcal{S}(\R ^3,\C )$ for all
$\mu$ and $\nu$. Let $\phi$ be its Lie transform. Denote $z'=\phi(z)$,
$z\equiv(\xi,f)$ and $z'\equiv(\xi ',f')$. Then there exist functions
$ G_{\mu,\nu} (z ), G_j (z)$ and a suffitiently small neighbourhood of
the orgin $\U^{-\kappa,-s}\subset
\Ph ^{-\kappa,-s}$, with the following three properties, which hold in
$\U^{-\kappa,-s}$.

\begin{itemize}
\item[1.]  $G_j, G_{\mu,\nu}\in C^\infty
(\U^{-\kappa,-s } , \mathbb{C}) $.  Actually such functions are
  analytic, but this will not be needed.
\item[2.] The transformation $\phi$ has the following structure:
\begin{eqnarray}
\label{lie.11.a}
\xi_j'&=&\xi_j+G_j(z) \\
\label{lie.11.b} f'&=& f+\sum_{\mu,\nu}G_{\mu,\nu}(z) \Psi_{\mu,\nu}
.
\end{eqnarray}

\item[3.] There are constants $C_{\tau,k,s}$   such that
\begin{equation}
\label{lie.11}
 \norma{z-\Tra(z)}_{\Ph ^{\kappa ,\tau}}\leq
 C_{\tau,k,s}  |\xi | ^{M_0}
 ( |\xi |+ \norma{f }
_{H^{-\kappa ,-s}} ).
\end{equation} Furthermore there are constants $c_{\kappa,\tau,k s}$
such that
\begin{eqnarray}
\label{lie.11.c}&
 |G_j(\xi,f)| \leq c_{\kappa,\tau,k ,s} |\xi | ^{M_0} ( |\xi |+
 \norma{f } _{H^{-\kappa ,-s}} ), \\ & \label{lie.11.d} |G_{\mu,\nu}
 (\xi,f)|\leq c_{\kappa ,s} |\xi | ^{{M_0}+1}.
\end{eqnarray}
\end{itemize}
\end{lemma}
\proof Recall $\phi =\phi ^1$. We set $z(t)=\phi ^t(z)=(\xi (t) , f(t)).$
The Hamilton equations of $\chi$  have the structure
\begin{equation}
\label{lie.4} \dot f= -\im \sum_{\mu,\nu}\xi^\mu\bar\xi^\nu
\Psi_{\mu,\nu}\ ,\quad \dot \xi_j=P_j(\xi)+\sum_{\mu,\nu}\tilde
P_{\mu,\nu}(\xi)\int_{\R^3}\Phib_{\mu,\nu}\cdot \beffe dx
\end{equation}
with suitable polynomials $P_j(\xi)$ homogeneous of degree ${M_0}+1$
and $\tilde P_{\mu,\nu}(\xi)$ homogeneous of degree ${M_0}
$.   By the existence and uniqueness theorem
for differential equations the solution exists up to time 1,
provided that the initial data are small enough. We consider
\eqref{lie.11}. For $t\in [0,1]$ we have for $\Ph $ equal to either
$\Ph ^{\kappa ,\tau}$ or  $\Ph ^{-\kappa ,-\tau}$
\begin{equation}
\label{lie.12bis} \begin{aligned} &  \| z(t)-z \| _{\Ph
}=\norma{\int_0^t X_{\chi}(z(t'))dt'}_{\Ph  }
\\& \le
 \widetilde{c}_{\kappa ,s} \sup _{0\le t'\le t} |\xi (t') | ^{M_0} ( |\xi (t')|+
 \norma{f(t') } _{H^{-\kappa ,-s}} ) .
\end{aligned}
\end{equation}
Then \eqref{lie.12bis} implies $|\xi (t )|+\norma{ f(t)
}_{H^{-\kappa ,-s} } \approx |\xi  |+\norma{ f  }_{H^{-\kappa ,-s} }
$ and $ |\xi (t )| \approx |\xi  | .$ Taking $t=1$ in the
rhs of \eqref{lie.12bis} we get
\begin{equation}
\label{lie.12} \begin{aligned} &  \| \phi (z) -z \| _{\Ph
}=\norma{\int_0^1 X_{\chi}(z(t'))dt'}_{\Ph  }
\\& \le
 c_{\kappa ,s}   |\xi   | ^{M_0} ( |\xi  |+
 \norma{f  } _{H^{-\kappa ,-s}} ) .
\end{aligned}
\end{equation}
\eqref{lie.12} is \eqref{lie.11}.
Any map $(\xi , f ) \to \xi '$ can be written in the form
 \eqref{lie.11.a}.    From the first of eq.\eqref{lie.4}, equation
\eqref{lie.11.b} holds with
$$ G_{\mu,\nu}(\xi (0),f (0)):=-\im
\int_0^1\xi^\mu(s,\xi(0),f(0))\bar \xi^\nu
(s,\xi(0),f(0))\Phi_{\mu\nu} ds\ .
$$ The $G_j$ in \eqref{lie.11.a} and the $G_{\mu,\nu}$ in \eqref{lie.11.b}
are analytic by the analyticity of flow $\phi^t(\xi,f )$, which is a
consequence of the analyticity of $X_{\chi }$ as a function defined
in $\Ph ^{-\kappa ,-s}$. \qed

\begin{lemma}
\label{kreg} Let $K\in C^k(\U^{1/2,0},\C)$, $k\geq 3$ satisfy
$|K(z)|\leq C\norma{z}^{M_1}$, and $\norma{dK(z)}_{\Ph^{-1/2,0}}\leq
C_1\norma{z}^{M_1-1}$, with $M_1\geq 2$. Let $\phi$ be the Lie
transform generated by the function $\chi$   of Lemma
\ref{lie_trans}. Then $K\circ \phi\in C^k(\Ph^{1/2,0},\R)$ and
$\{K,\chi\}\in C^{k-1}(\U^{1/2,0},\R)$. Furthermore one has
\begin{eqnarray}
\label{kreg1}
\left|K(\phi(z))\right|\leq
C\norma{z}^{M_1}\ ,
\\
\label{kreg2} \left|K(\phi(z))-K(z)\right|\leq C\norma{z}^{M_0+M_1}.
\end{eqnarray}
\end{lemma}
\proof \eqref{kreg1} is an elementary consequence of \eqref{lie.11}. We have
\begin{equation} \begin{aligned}& \left|K(\phi(z))-K(z)\right|\leq
    \|    \phi(z)-z \|  _{\Ph^{1/2,0}}  \sup _{t\in [0,1]} \| d K
    (z +t(\phi (z) -z))\| _{\Ph^{-1/2,0}}
\\& \le C\norma{z}^{M_o+M_1} ,
  \end{aligned} \nonumber
\end{equation}
by  $\| d K (z  )\| _{\Ph^{-1/2,0}} \le  C_1\norma{z}^{M_1-1}  $ and by \eqref{lie.11}.
   \qed

The next lemma is elementary.
\begin{lemma}
\label{smo.meno} Let $K\in C^{\infty}(\U^{-k,-s},\C )$, where
$\U^{-k,-s}\subset \Ph^{-k,-s}$, with some $s\geq 0$, $k\geq 0$.
Then one has $X_{K}\in C^{\infty}(\U^{-k,-s},\Ph^{k,s}) $.
\end{lemma}

\subsection{Normal form}

\begin{definition}
\label{d.1} A polynomial $Z$ is in normal form if  we have
\begin{equation}
\label{e.12} Z=Z_0+Z_1
\end{equation}
where:   $Z_1$ is a linear combination of monomials of the form
\begin{equation}
\label{e.12a} \xi^\mu\bar\xi^\nu\int \Phi(x)f(x) dx\ ,\quad
\xi^{\mu'}\bar\xi^{\nu'}\int \Phi(x)\bar f (x)dx
\end{equation}
with   indexes  satisfying
\begin{equation}
\label{e.12b} \omega\cdot(\mu-\nu)<-m \ ,\quad \omega\cdot (
\mu'-\nu ')>m \ ,
\end{equation}
and $\Phi\in  \mathcal{S}(\R ^3, \C )$; $Z_0 $ is independent of
 $f$ and is a linear combination of monomials
 $\xi^\mu\bar\xi^{\nu}$ satisfying
\begin{equation}
\label{e.13} \left\{H_L,\xi^\mu\bar\xi^{\nu}\right\}=0.
\end{equation}
\end{definition}

\begin{remark} Equation \eqref{e.13} is equivalent to $\omega \cdot
  (\mu -\nu )=0$, see Lemma \ref{homdia} below.
\end{remark}
\begin{remark}\label{alfabeta} By (H5),
$ \omega\cdot(\mu-\nu)=0$ implies
$\left|\mu\right|=\left|\nu\right|$ if $|\mu +\nu|\le 2N_1+3$.
\end{remark}

\begin{theorem}
\label{main} For any  $k>0$ and  $s>0$ and for any    integer $r $ with $0\le
r \le 2N$ there exist open neighborhoods of the origin
$\U_{r,k,s}\subset\Ph^{1/2,0}$, and $\U^{-k,-s}_r\subset \Ph^{-k,-s}$, and
an analytic canonical transformation $\Tr_r:\U_{r,k,s}\to\Ph^{1/2,0}$ with the following properties. First of all $\Tr_r$   does not depend on $(k ,s )$
 in the sense that,
 given another pair $(k',s')$, the transformations coincide in   $\U_{r,k,s}\cap \U_{r,k',s'}$. Secondly,
$\Tr_r$ puts the system in normal form up to order $r+4$, namely we have
\begin{equation}
\label{eq:bir1} H^{(r)}:=H\circ \Tr_r=H_L+Z^{(r)}+\resto^{(r)}
\end{equation}
where:
\begin{itemize}
\item[(i)]$Z^{(r)}$ is a polynomial of degree $r+3$ which is in normal
  form; furthermore, when we expand
\begin{equation}
\label{ze.1}
Z^{(r)}_1(\xi,f)=\sum_{\mu,\nu}\xi^{\mu}\bar\xi^\nu
\int_{\R^3}\Phi_{\mu\nu}f dx+ \sum_{\mu,\nu}\bar \xi^{\mu}\xi^\nu
\int_{\R^3}\bar \Phi_{\mu\nu}\bar f dx
\end{equation}
we have, for $\beta_{|\mu|}:=\beta^{(|\mu|)}(0)$,
$\varphi^\mu=\prod_{j}\varphi_j^{\mu_j}$ and similarly $\omega^\mu
=\prod_{j}\omega _j^{\mu_j}$,
\begin{equation}
\label{ze.2}\Sc(\R^3,\C ) \ni \Phi_{\mu0}=\frac{2^{-\frac{|\mu
|}{2}}}{\mu !}\beta _{| \mu |+1 }
  \frac{B^{-\frac{1}{2}}(\varphi ^{\mu }) (x)} {\sqrt{\omega ^{\mu
  }}}+\tilde \Phi_{\mu0}
\end{equation}
with $\tilde \Phi_{\mu0}  = \tilde \Phi_{\mu0}
(m,\beta_4,...,\beta_{|\mu|}) $  piecewise smooth in
$(m,\beta_4,...,\beta_{|\mu|})$, with values in $\Sc(\R^3,\C) $; the functions $\Phi_{\mu\nu}(x)$ belong to  $\mathcal{S}(\R ^3, \C ^2); $

\item[(ii)]     $\Tr _r$ has the structure \eqref{lie.11.a},
  \eqref{lie.11.b},  $\uno-\Tr_r$ extends into an
  analytic map from $\U^{-k,-s}_r$ to $\Ph^{k,s}$ and
\begin{equation}
\label{def1} \norma{z-\Tr_r(z)} _{\Ph ^{k,s}}
 \leq C\norma{z}_{\Ph
^{-k,-s} } ^3 \,  ;
\end{equation}
\item[(iii)] we have $\resto^{(r)} = \sum _{d=0}^4\resto^{(r)}_d$
with the following properties:
\begin{itemize}
\item[(iii.0)] we have

\begin{equation} \resto^{(r)}_0=
\sum _{|\mu +\nu |  = r+4 } \xi ^\mu \overline{\xi}^\nu \int
_{\mathbb{R}^3}a_{\mu \nu }^{(r)}(x,z,\Re  B^{-\frac{1}{2}}f (x)  ) dx
\nonumber
\end{equation}
and   $a_{\mu \nu }^{(r)}$ is such that the map
\begin{equation}\label{eq:coeff a}
\U^{-k,-s}\times \R \ni(z,w)\mapsto a_{\mu \nu }^{(r)}(.,z, w  )\in
H^{k,s} \quad\text{is}\ C^\infty
\end{equation}
\item[(iii.1)] we have
\begin{equation} \resto^{(r)}_1=
\sum _{|\mu +\nu |  = r +3 }\xi^\mu \overline{\xi}^\nu \int
_{\mathbb{R}^3}   {\bf \Lambda }  _{\mu \nu
}^{(r)}(x,z,\Re B^{-\frac{1}{2}}f(x) ) \cdot B^{-\frac{1}{2}}\beffe(x)
dx \nonumber
\end{equation}
where  the map
\begin{equation}\label{eq:coeff beffe}
\U^{-k,-s}\times \R \ni(z,w)\mapsto {\bf \Lambda }  _{\mu \nu
}^{(r)}(\cdot,z, w  )\in (H^{k,s})^2 \quad\text{is}\ C^\infty
\end{equation}
\item[(iii.2-3)]
for $d=2,3$,  we have
\begin{equation}
\label{rd03}
\resto^{(r)}_d=\int_{\R^3}F^{(r)}_d(x,z,\Re B^{-\frac{1}{2}}f(x))[U(x)]^ddx\
,
\end{equation}
where $U=B^{-1/2}(f+\bar f)$ where the map
\begin{equation}\label{eq:coeff F1}
\U^{-k,-s}\times \R \ni(z,w)\mapsto F^{(r)}_d(\cdot,z, w  )\in
H^{k,s} (\R ^3, \C )\quad\text{is}\ C^\infty
\end{equation}
and furthermore  we have
\begin{equation}\label{eq:coeff F}   \|  F_2^{(r)} (\cdot ,z ,w )\| _{  H^{k,s} (\R ^3, \C )  } \le C   |\xi   |   ;
\end{equation}

\item[(iii.4)] for $d=4$ we have
\begin{equation}
\label{rd4} \resto^{(r)}_4=\int_{\R^3}F _4( x,\Tr _r(z))[U(x)]^4dx\
,
\end{equation}
where $F _4 (x,z)= F_4(x,\xi ,U)$ is the function in
\eqref{taylor2}.
\end{itemize}
\end{itemize}
\end{theorem}


\subsection{The Homological Equation}\label{homs}

Let $K(\xi,\bar\xi,f,\bar f)$ be a homogeneous polynomial of degree
$M_1$ of the form
\begin{eqnarray}
\label{eq.stru.1} K =\sum_{|\mu|+|\nu|=M_1}K_{\mu\nu}
\xi^\mu\bar\xi^\nu+ \sum_{|\mu'|+|\nu'|=M_1-1}
\xi^{\mu'}\bar\xi^{\nu'}\int \Phi_{\mu'\nu'} f
\\
\nonumber
+
\sum_{|\mu''|+|\nu''|=M_1-1} \xi^{\mu''}\bar\xi^{\nu''}\int
\Psi_{\mu''\nu'' } \bar f\ ,
\end{eqnarray}
with functions $ \Phi_{\mu'\nu'}, \Psi_{\mu''\nu''}\in
\mathcal{S}(\R ^3, \C)$. A key step in the proof of Theorem
\ref{main} consists in solving  (i.e. finding $\chi$ and $Z$) with
$Z$   in normal form, the homological equation
\begin{equation}
\label{homo} \left\{   H_L , \chi \right\}+Z=K\ .
\end{equation}
  To solve \eqref{homo} we first define $Z$ to be the r.h.s. of \eqref{eq.stru.1}
restricting the sum to the indexes such that
\begin{equation}
\label{e.12ba} \omega\cdot(\mu-\nu)=0\ ,\quad \omega\cdot
(\nu'-\mu')>m \ ,\quad \omega\cdot (\mu''- \nu'')> m\ ,
\end{equation}
i.e. the indexes of  the  normal form condition. We introduce the {\it
homological operator}
\begin{equation}
\label{homop} \lie\chi:=\left\{   H_L , \chi \right\}
\end{equation}

\begin{lemma}
\label{homdia} We have: \begin{eqnarray} \label{homo.a.1} \lie (
\xi^\mu\bar\xi^\nu ) = - \im \omega\cdot(\mu-\nu
)\xi^\mu\bar\xi^\nu\ ,
\\
\label{homo.a.2} \lie(\xi^\mu\bar\xi^\nu\int \Phi f)= -\im
\xi^\mu\bar\xi^\nu\int  f ( B- \omega\cdot (\nu-\mu)   )\Phi  \ ,
\\ \label{homo.a.3}
\lie(\xi^\mu\bar\xi^\nu\int \Phi \bar f)=  \im
\xi^\mu\bar\xi^\nu\int \bar f ( B - \omega\cdot (\mu-\nu)   )\Phi  \
.
\end{eqnarray}
\end{lemma}
\proof Indeed, using \eqref{e.9b}, \eqref{homo.a.1} follows by
\begin{eqnarray*}
\lie (\xi^\mu\bar\xi^\nu ) = \im\sum_j \omega_j\left( \bar
\xi_j\frac{\partial}{\partial \bar \xi_j}-
 {\xi}_j\frac{\partial}{\partial {\xi}_j}\right)
\xi^\mu\bar\xi^\nu = \im \omega\cdot(\nu-\mu)\xi^\mu\bar\xi^\nu \ .
\end{eqnarray*}   \eqref{homo.a.2}--\eqref{homo.a.3} follow from
\eqref{e.9b}, \eqref{homo.a.1},    $\lie(\langle  \Phi
,\overline{f}\rangle ) = \im \langle  \Phi , B\overline{f}\rangle$,
$\lie(\langle  \Phi , {f}\rangle ) = -\im \langle  \Phi , B
{f}\rangle$ and selfadjointness of $B$. \qed

\medskip

 For $\omega\cdot (\mu-\nu)<m$ we set
\begin{equation}
\label{res} R_{\mu\nu}:=(B-\omega\cdot (\mu-\nu) )^{-1} \ .
\end{equation}
Notice   that $(B-\lambda)^{-1}$ is a real operator for $\lambda<m$.
Then, Lemma \ref{homdia} yields immediately:

\begin{lemma}
\label{sol.homo} Let $K$ be a polynomial  as in \eqref{eq.stru.1};
define $Z$ as above and $\chi :=$
\begin{equation}
\label{solhomo}    \sum_{\alpha ,\beta}\frac{\im K_{\alpha \beta}
\xi^  {\alpha   }\bar\xi ^ {\beta  }}{
  \omega\cdot(\alpha -\beta )}  +
   \im  \sum_{\mu  ,\nu  }
 \xi^{\mu }\bar\xi^ {\nu   }\int  fR_{\nu  \mu}  \Phi_{\mu  \nu}-
 \im \sum_{\mu ' ,\nu ' }
  \xi^{\mu ' }\bar\xi^{\nu ' }\int
 \bar f R_{\mu '  \nu' }\Psi _{\mu' \nu ' }
\end{equation}
with the sum restricted  to indexes of the sum \eqref{eq.stru.1}
such that
\begin{equation}
\label{e.12bac} \omega\cdot(\alpha -\beta )\not=0\ ,\quad
\omega\cdot (\nu -\mu )<m \ ,\quad \omega\cdot (\mu' - \nu' )< m\ .
\end{equation}
Then equality \eqref{homo}  is true for this choice of $\chi$ and
$Z$. Furthermore, if $K_{\mu \nu}= \overline{K}_{\nu \mu}$ and $\Psi
 _{\mu \nu}= \overline{\Phi }_{\nu \mu}$, also the coefficients in
\eqref{solhomo} and in the sum defining $Z$ satisfy this property.
\end{lemma}
We also need the following regularity result, proved in Appendix \ref{app.resolv}
 at the end of the paper.
\begin{lemma}
\label{sol.homo1} Suppose (H1)--(H2),    $\Phi =P_c\Phi$ and $\Phi
\in \mathcal{S}(\R ^3, \C )$. Then:
\begin{itemize}
\item[(1)] for
 $\lambda <m$  we have $
(B-\lambda )^{-1} \Phi \in \mathcal{S}(\R ^3, \C )$; \item[(2)] for
any $l\in \R $ we have $B^l \Phi \in \mathcal{S}(\R ^3, \C )$.
\end{itemize}
\end{lemma}

\subsection{ Proof of Theorem \ref{main}}

\noindent {\bf Proof of Theorem \ref{main}.} By Lemma
\ref{lemma:H_P}, $H$ satisfies   assumptions and conclusions of
Theorem \ref{main} with $r=0$, $\Tr_0\equiv\uno$,
$\resto^{(0)}:=H_P$, $Z^{(0)}=0$. We now assume that the theorem is
true for $r$ and   prove it for $r+1$. Define
\begin{eqnarray}
\label{r00} &\resto^{(r)}_{02} = \resto^{(r)}_{0}-\sum _{|\mu +\nu
|  = r+4 } \xi^\mu \overline{\xi}^\nu  \int _{\mathbb{R}^3}a_{\mu
\nu }^{(r)}(x,0 ,0 ) dx,
\\
\label{r11} & \resto^{(r)}_{12} =\resto^{(r)}_{1}-\sum _{|\mu +\nu
| = r +3 }\xi^\mu \overline{\xi}^\nu \int _{\mathbb{R}^3} \Phib
_{\mu \nu }^{(r)}(x   )\cdot \beffe (x) dx,
\end{eqnarray}
with  $\Phib _{\mu \nu }^{(r)}(x   )=B^{-\frac{1}{2}} {\bf \Lambda}
_{\mu \nu }^{(r)}(x,0 ,0 ) $. Notice that even though the rhs of
\eqref{eq:bir1} can depend on the pair $(k,s)$,
 the terms      $ {\bf \Lambda} _{\mu
\nu }^{(r)}(x,0 ,0  )=\frac{1}{\mu !\nu !}\partial ^{\mu}_{
{\xi}}\partial ^{\nu}_{\overline{\xi}}\nabla _f H^{(r)}(0)$ are
independent of $(k,s)$ (because of the independence on $(k,s)$ of $\Tr _r$, and hence of $H^{(r)}$, as a germ at the origin). Hence $ {\bf \Lambda} _{\mu \nu }^{(r)}(x,0
,0 )\in \mathcal{S}(\R ^3, \C ^2)$. Then  $\Phib _{\mu \nu
}^{(r)}(x   ) \in \mathcal{S}(\R ^3,\C )$ by  Lemma
\ref{sol.homo1}. We have
 \begin{equation} \label{eq:r00} \begin{aligned} &
 \resto^{(r)}_{02} + \resto^{(r)}_{12} =\sum _{|\mu +\nu
|  = r+5 } \xi^\mu \overline{\xi}^\nu  \int
_{\mathbb{R}^3}\widetilde{a}_{\mu \nu }^{(r)}(x,z ,0  ) dx+\\& \sum
_{|\mu +\nu |  = r+4  }\xi^\mu \overline{\xi}^\nu \int
_{\mathbb{R}^3}  \widetilde{{\bf \Lambda}}_{\mu \nu }^{(r)}(x,z
,\Re B^{-\frac{1}{2}} f(x) ) \cdot B^{-\frac{1}{2}} \beffe(x) dx+ \\&
\sum _{|\mu +\nu | = r +3  }\xi ^\mu \overline{\xi}^\nu \int
_{\mathbb{R}^3} \widetilde{F}_{2\mu\nu }^{(r)}(x,z
,\Re B^{-\frac{1}{2}} f(x)  )
 \cdot  \left (B^{-\frac{1}{2}} \beffe (x)\right ) ^{  2}  dx,
 \end{aligned}
  \end{equation}
with $\widetilde{a}_{\mu \nu }^{(r)}$ satisfying \eqref{eq:coeff
a}, $\widetilde{{\bf \Lambda}}_{\mu \nu }^{(r)}$
\eqref{eq:coeff beffe} and with $\widetilde{F}_{2\mu\nu }^{(r)}$
such that the map
\begin{equation}  \begin{aligned} &
\U^{-k,-s}\times \R \ni(z,w)\mapsto \tilde F_{2\mu\nu}^{r}(.,z, w  )\in
H^{k,s} \quad\text{is}\ C^\infty
\end{aligned}\nonumber
\end{equation}
Set
\begin{equation}
\label{pr.3}{K}_{r+1} := \sum _{|\mu +\nu |  = r+4 } \xi^\mu
\overline{\xi}^\nu \int _{\mathbb{R}^3}a_{\mu \nu }^{(r)}(x,0 ,0 )
dx+ \sum _{|\mu +\nu | = r +3 }\xi^\mu \overline{\xi}^\nu \int
_{\mathbb{R}^3} \Phib _{\mu \nu }^{(r)}(x  )\cdot \beffe (x)
dx.\nonumber
\end{equation}
${K}_{r+1} $ is real valued, so in particular its coefficients satisfy    the last sentence of Lemma \ref{sol.homo}.
We can apply Lemma \ref{sol.homo} and denote by $\chi_{r+1}$ and
$Z_{r+1}$ the solutions of the homological equation
$$
\left\{    H_L ,  \chi_{r+1}  \right\}+Z_{r+1}=K_{r+1}\ .
$$ Let $\phi_{r+1}$ be the Lie transform generated by
$\chi_{r+1}$. The discussion in Remark \ref{rem.re} holds.
Let $\U_{r+1}$, $\U_{r+1}^{-k,-s}$ be such that
$\phi_{r+1}(\U_{r+1})\subset \U_{r}$ and
$\phi_{r+1}(\U_{r+1}^{-k,-s})\subset \U_{r}^{-k,-s}$. Denote $(\xi
',f')=\phi_{r+1}(\xi ,f)$. Then
$f'=f+\sum_{\mu\nu}\Psi^{(r+1)}_{\mu\nu} G^{(r+1)}_{\mu\nu}(z)$, with
$G^{(r+1)}_{\mu,\nu}$ described by Lemma \ref{lie_trans} and
$\Psi^{(r+1)}_{\mu\nu}\in\Sc(\R^3,\C)$.
Denote
$$G_U:=B^{-1/2}\sum_{\mu\nu}(\Psi^{(r+1)}_{\mu\nu} G^{(r+1)}_{\mu\nu} +\bar
\Psi^{(r+1)}_{\mu\nu} \bar G^{(r+1)}_{\mu\nu})\ . $$ Recall \eqref{lie.11.b}
and \eqref{lie.11.d},
which imply
\begin{equation}
\label{dim.12}
\norma{G_U(z)}_{H^{k,s}}\leq C|\xi |^{r+3} \ .
\end{equation}
We define by induction    $\Tr _0= \uno$, $\Tr  _{r+1}= \Tr _{r
}\circ\Tra  _{r+1}$. Then  \eqref{lie.11} implies claim (ii).

We will now   prove that
$$
H^{(r+1)}:=H^{(r)}\circ\phi_{r+1}\equiv H\circ(\Tr_r\circ\phi_{r+1})
\equiv H\circ \Tr_{r+1} \ ,
$$
has the desired structure. Write
\begin{eqnarray}
\label{n1.1} H^{(r)}\circ\Tra_{r+1}&=& H_L+Z^{(r)} + Z_{r+1}
\\
\label{n1.3}
&+& (Z^{(r)}\circ\Tra_{r+1}-Z^{(r)})\\
\label{n1.6} &+& K_{r+1}\circ\Tra_{r+1}-K_{r+1}
\\
\label{n1.4} &+& H_L\circ\Tra_{r+1}- \left(H_L+\left\{ \chi_{r+1},
 H_L  \right\}\right)
\\
\label{n1.5} &+&
(\resto^{(r)}_{02}+\resto^{(r)}_{12})\circ\Tra_{r+1}
\\
\label{n1.7} &+&   \resto ^{(r)}_2\circ\Tra_{r+1}
\\
\label{n1.8} &+&   \resto
^{(r)}_3\circ\Tra_{r+1}  \\
\label{n1.9} &+& \resto ^{(r)}_4\circ\Tra_{r+1} \ .
\end{eqnarray}
  $Z^{(r+1)}:=Z^{(r)}+Z_{r+1}$ is in normal form and of the desired
  degree.
    We   study  now
\eqref{n1.7} and \eqref{n1.8}.  For $d=2,3$, expanding $(U+G_U)^d$  one has
\begin{equation}
 \begin{aligned}
  &\left(\resto ^{(r)}_d\circ\Tra_{r+1}
\right)(z)=\sum_{j=0}^{d}\left(
\begin{matrix}
d\\ j
\end{matrix}
\right) \int_{\R^3} F^{r}_d\left(... \right)
[G_U(z)]^{d-j}[U(x)]^j=: \sum_{j=0}^{ d}H_{dj}\ ,\\& \text{where } F^{r}_d\left(... \right)= F^{r}_d\left(x,\phi_{r+1}(z), \Re  B ^{-\frac{1}{2}} (f+
\sum_{\mu\nu}\Psi_{\mu\nu}^{(r+1)}G^{(r+1)}_{\mu\nu}(z) )(x) \right) .
 \end{aligned}\nonumber
\end{equation}
Each of the functions $H_{dj}$ has the structure (iii.0-iii.4).
Similarly
\begin{eqnarray}
\label{pe.2} (\resto ^{(r)}_4\circ\Tra_{r+1})(z)= \sum_{d=0}^{4}\left(
\begin{matrix}
4 \\ d
\end{matrix}
\right)  \int_{\R^3}   F _4(x,\phi_{r+1}(z)  )
[G_U(z)]^{4-d}[U(x)]^{d }.\nonumber
\end{eqnarray}
Each term  with $d\le 3$ can be absorbed in  $\resto^{(r+1)}_d$.
 For $d=4$ we get
(iii.4).  \eqref{n1.5} can be treated similarly. Notice   that, by \eqref{dim.12}, all the contributions
to $\resto ^{(r+1)}_2$  from the  $H_{dj}$ satisfy \eqref{eq:coeff F}. The same is true
for the contributions coming from \eqref{n1.5}, i.e. from the last line of \eqref{eq:r00}, and from \eqref{n1.9}.

 By $K_{r+1}\in
C^{\infty}(\U_{r}^{-k,-s})$, the term \eqref{n1.6} can be included
in $\resto^{(r+1)}_0$, with the vanishing properties  at $\xi =0$
and $f=0$  guaranteed by \eqref{kreg2}.   \eqref{n1.3} can be
treated exactly in the same way. We prove now that \eqref{n1.4} can
be included in $\resto ^{(r+1)}_0 $. We write
\begin{eqnarray}
 & H_L\circ\phi_{r+1}-(H_L+\left\{ \chi_{r+1},
 H_L  \right\})=
\int_0^1 \frac{t^2}{2!} \frac{d^2}{dt^2}\left(
H\circ\phi_{r+1}^t\right) dt \label{hllis}
\\&
 \nonumber
=\int_0^1 \frac{t^2}{2!} \left\{
 \chi_{r+1} ,\left\{ \chi_{r+1},
 H_L  \right\} \right\}\circ\phi_{r+1}^t dt= \int_0^1
\frac{t^2}{2!}\left\{
 \chi_{r+1}, Z_{r+1}-K_{r+1}  \right\}\circ\phi_{r+1}^t
dt.
\end{eqnarray}
This shows that \eqref{n1.3} is in $C^{\infty}(\U_{r+1}^{-k,-s})$,
with   vanishing properties  at $z=0$ which allow to absorb it in
$\resto ^{(r+1)}_0.$

We prove   equation \eqref{ze.2}. Consider $\Phi_{\mu 0}$ with
$|\mu| =r+2$. Then
\begin{equation} \label{eq:phimu0}\mu !\Phi_{\mu 0}    = \partial ^\mu _\xi \nabla  _f
 H^{(r)}(0)  .\nonumber
\end{equation}
We have
\begin{eqnarray}&  \partial ^\mu _\xi \nabla  _f
 H^{(r )}(0) =  \partial ^\mu _\xi \nabla   _f
 H^{(0)}(0)+  \partial ^\mu _\xi \nabla   _f
  \left [  H^{(0)}\circ  \Tr _{r }- H^{(0)} \right ](0)  =\nonumber \\&
  \label{m.35}2^{\frac{r +3}{2}}\beta ^{(r +4) } (0)
  \frac{B^{-\frac{1}{2}}(\varphi  ^{\mu }) (x)} {\sqrt{\omega  ^{\mu }}} +
    \partial ^\mu _\xi \nabla   _f
  \left [H^{(0)}\circ  \Tr _{r }- H^{(0)} \right ](0)
\end{eqnarray}
where the first term in the right hand side is obtained by Lemma
\ref{lemma:H_P}. So we need to show that the last term in
\eqref{m.35} is like the reminder in \eqref{ze.2}. First of all
notice that if we consider the embedding $I_k:\Ph
^{k,0}\hookrightarrow \Ph ^{ \frac{1}{2},0}$ for $k>1/2$ with
$I_k(z)=z$, we have $\partial ^\mu _\xi \nabla _f H^{(r)}(0)
=\partial ^\mu _\xi \nabla _f [ H^{(r) }\circ I_k](0) $ for any
$\mu$. In other words, it is enough that we prove our formula
restricting the Hamiltonians on $\Ph ^{k,0}$ for $k$ large. We prove
that $d^{r+4}\left [H^{(0)}\circ \Tr _{r }- H^{(0)} \right ](0)$ is
a smooth function of $ ( m, \beta _4, ...,\beta_{r +3} )$, where
$\beta _l:=\beta ^{(l)}(0)$, with $m>\lambda _1$ such that
(H3)--(H5) are satisfied. We can apply the chain rule and obtain the
standard formula
\begin{eqnarray}\label{chainrule}& d^{r+4}  (H ^{(0)}\circ
 \Tr  _r)(0) =\sum _\alpha c_\alpha
 (d^{|\alpha | }H ^{(0)})(0)  \left (\otimes _{j=1}^{r+4} ( d^j \Tr
_r(0))^{\alpha _j}\right )
\end{eqnarray}
with $\sum _{j=1}^{r+4} j\alpha _j=r+4$ and $c_\alpha$ appropriate
universal constants. Insert the decomposition $ \Tr _{r }= \uno +
\widetilde{\Tr} _{r }$ into \eqref{chainrule}. Then $d^{r+4}  (H
^{(0)}\circ
 \Tr  _r)(0) =d^{r+4}   H ^{(0)}+ \mathcal{E}$ where  $\mathcal{E}$
 is a sum of terms of the form
\begin{eqnarray}\label{chainrule1}& c_\alpha
 (d^{|\alpha | }H ^{(0)})(0) \left ( \uno ^{ \otimes \alpha_0}\otimes _{j=1}^{r+4}
  ( d^j \widetilde{\Tr}
_r(0))^{\widetilde{\alpha} _j} \right )
\end{eqnarray}
with at least one $\widetilde{\alpha} _j>0 $ and for some $\alpha
_0\ge 0$.  By $d ^j \widetilde{\Tr} _r(0)=0$ for $0\le j \le 2$ we
have $\widetilde{\alpha} _1=\widetilde{\alpha} _2=0$ and so $\alpha
_j=\widetilde{\alpha }_j>0$ for some $j\ge 3$. Hence the terms in
\eqref{chainrule1} are such that $|\alpha |< r+4 .$ $(d^{j }H
^{(0)})(0) $ for $j<r+4$ is a smooth function of $( m, \beta _4,
...,\beta_{r +3} )$. Indeed, if we reverse the change of variables
\eqref{e.8}, $d^{j } H ^{(0)} (0)= \beta _{j} $ for all $j $. By
induction it is elementary to show that $\widetilde{\Tr} _{r }(z)=
\widetilde{\Tr} _{r }(z,m, \beta _4, ...,\beta_{r +3} )$ is a
smooth
function of all its arguments.  In particular it is smooth also in
  $m$ for all values such that  $m>\lambda _1$ and that
(H3)--(H5) are satisfied.  Indeed $\widetilde{\Tr} _{0}\equiv 0$,
$\widetilde{\Tr} _{r }$ depends on the vector field $K_{r }$ which
in turn is a smooth function of
\begin{equation}\label{chainrule2}
\partial ^{\nu}_\zeta \nabla _\beffe ^{j} H ^{(r
-1)} (0) \text{ with $|\nu |+j=r +3$ and $j\le 1$.}\end{equation} By
 induction, \eqref{chainrule2} is a smooth function of $( m, \beta _4,
 ...,\beta_{r +3} )$ with $m>\lambda _1$ such that
(H3)--(H5) are satisfied. Hence we have also proved property (i) of
 Theorem \ref{main}. \qed

\section{Dynamics of the normal form}
\label{model}

Before giving the proof of Theorem \ref{theorem-1.1} we outline the
main features of the dynamics generated by the normalized system and
we discuss the nondegeneracy assumption. Our main idea has been to
normalize through canonical transformations. Hence we have preserved
the Hamiltonian nature of the system. We now proceed exactly as in the
literature, with the difference that at the end we can show the
positive semidefiniteness of some key coefficients, see Lemma
\ref{l.dis}. This semidefiniteness is in the literature either proved
in the special case $N=1$, or in very special cases.

  In the sequel we assume
that the time $t $ is positive.  Due to the time reversal invariance of
the equations,  this is not restrictive. We
consider $r=2N$. We neglect $\resto^{(2N)}$ and consider the
Hamiltonian
\begin{equation}
\label{m.1} H_{nf}:=H_L(\xi,\beffe)+Z_0(\xi)+Z_1(\xi,\beffe)\
.
\end{equation}
We show later that the addition of $\resto^{(2N)}$ to $H_{nf}$ does
not change the qualitative  features of  the dynamics of  the
simplified system considered in this section.  $Z_0$ and $Z_1$ are as
in Definition \ref{d.1}, where
\begin{eqnarray}
\label{m.2} Z_1(\xi,\bar\xi,f,\bar f):= \langle G,f\rangle+\langle
\bar G,\bar f\rangle\ , \\ G:=\sum_{\mu , \nu}\xi^\mu\bar \xi^\nu
\Phi_{\mu\nu}\ ,\quad \bar G=\overline{\sum_{\mu, \nu }\xi^\mu\bar
  \xi^\nu \Phi_{\mu\nu}},
\end{eqnarray}
$\Phi_{\mu\nu}\in  \mathcal{S}(\R ^3, \C)$, with $\mu,\nu$ such that
\begin{equation}
\label{m.3} 2\leq |\mu|+|\nu|\leq 2N+2\ ,\quad \omega\cdot
(\mu-\nu)<-m \ .
\end{equation}
The Hamilton equations of this system are given by
\begin{eqnarray}
\label{m.4} \dot f&=&-\im (Bf+\bar G)\ ,
\\
\label{m.5} \dot \xi_k&=& -\im  \omega_k \xi_k- \im \frac{\partial
Z_0}{\partial
  \bar \xi_k}-\im \left\langle \frac{\partial G}{\partial
  \bar \xi_k},f\right\rangle-\im \left\langle \frac{\partial \bar G}{\partial
  \bar \xi_k},\bar f\right\rangle
\end{eqnarray}
We prove later that $f$ is asymptotically free  in the dynamics of
the full system.
We need to examine
  in detail $f$ in order to extract its main contribution to the
equations for the $\xi _k$.  Hence we decouple further the dynamics
of the discrete modes and the continuous ones, following the
literature, see for instance \cite{cuccagnamizumachi} and references
therein. We do not change coordinates as in  the previous procedure,
since by the resonance between continuous and discrete spectrum the
Hamiltonian is not well defined in terms of the new decoupled
variables. So, as in the literature, we work at the level of vector
fields and look for a function $Y=Y(\xi,\bar \xi)$ such that the new
variable
\begin{equation}
\label{m.6}
g:=f+\bar Y
\end{equation}
is decoupled up to higher order terms from the discrete variables.
Substitution in   equation \eqref{m.4} yields
\begin{eqnarray}
\label{m.7.a} \dot g=-\im Bg -\im \left\{\bar G - \left[ B-\sum_k
  \left(\omega_k\xi_k\frac{\partial }{\partial \xi_k}-\omega_k\bar
  \xi_k\frac{\partial }{\bar \partial \xi_k} \right) \right] \bar Y
\right\}+{\rm h.o.t.}
\end{eqnarray}
where h.o.t. denotes terms which are either at least linear in
$\beffe$ or of sufficiently high degree   in $\xi$ (that is,
monomials $\xi ^\mu \bar \xi ^\nu$ with    $|\mu +\nu |>2N+2$).
  We want $Y$ such that the curly bracket vanishes. Write
\begin{equation}
\label{m.8} \bar Y:=\sum_{\substack{2\leq |\mu|+|\nu|\leq 2N+3\\
\omega\cdot (\mu-\nu)>m}}\bar Y_{\mu\nu}(x)\xi^\mu\bar \xi^{\nu}.
\end{equation}
The  vanishing of the curly bracket  in \eqref{m.7.a} is equivalent
to
\begin{equation}
\label{m.9} (B-\omega\cdot(\mu-\nu))\bar Y_{\mu\nu}=\bar \Phi_{ \nu
\mu }\, .
\end{equation}
Since $\omega\cdot(\mu-\nu)\in\sigma(B)$ we have to regularize the
resolvent. We set
\begin{equation}
\label{m.10}
R^{\pm}_{\mu\nu}:=\lim_{\epsilon\to0^+}(B-(\mu-\nu)\cdot \omega\mp
\im \epsilon )^{-1}.
\end{equation}
Now, in the sequel it is important that $t\ge 0$. We define
\begin{equation}
\label{m.11} \bar Y_{\mu\nu} = R^+_{\mu\nu} \bar \Phi_{ \nu \mu }
\quad \text{and} \quad Y_{\mu\nu} =\overline{R^+_{\mu\nu} \bar
\Phi_{ \nu \mu }} =R^-_{\mu\nu} \Phi_{ \nu \mu }\ .
\end{equation}
 \begin{lemma}
\label{lem:m.11} We have $Y_{\mu\nu}\in L^{2,-s}$
for all $s>1/2$, and thus also  $g\in L^{2,-s}$
for all $s >1/2$. \end{lemma}
\proof Follows immediately from Lemma \ref{lemma-regularization}
in  Appendix \ref{app.resolv}.  \qed

We substitute \eqref{m.6} in the equations for $\xi$, namely
\eqref{m.5}. Then we get
\begin{eqnarray}
\label{m.12} \dot \xi_k&=& -\im  \omega_k \xi_k- \im  \frac{\partial
Z_0}{\partial
   \overline{\xi}_k}+\im  \left\langle \frac{\partial G}{\partial
    \overline{\xi}_k},\bar Y\right\rangle+\im
    \left\langle \frac{\partial \bar
  G}{\partial
    \overline{\xi}_k},Y \right\rangle
\\
 \label{m.12a} && -\im  \left\langle \frac{\partial G}{\partial
    \overline{ \xi}_k},g\right\rangle-\im  \left\langle \frac{\partial \bar G}{\partial
    \overline{\xi}_k},\bar g\right\rangle\ .
\end{eqnarray}
We   show in the next section that $g$ is negligible.  So we neglect
\eqref{m.12a}.    A simple explicit computation using \eqref{m.2},
\eqref{m.8} and \eqref{m.11},  shows that the system \eqref{m.12} is
of  the form
\begin{eqnarray}
\label{m.13} &\dot \xi_k =  -\im   \omega_k \xi_k- \im
\frac{\partial Z_0}{\partial
    \overline{\xi}_k}
\\
\label{m.14} &  +\im     \sum_{ \substack{ \omega \cdot (\nu - \mu
)>m \\ \omega \cdot (\mu '- \nu ')>m
 } } \frac{\xi^{\mu+\mu'}  \overline{\xi}^{\nu'+\nu}}
{ \overline{ \xi}_k} \nu_kc_{\mu\nu\mu'\nu'}+\\
\label{m.14bis} &  +\im    \sum_{ \substack{ \omega \cdot (\nu - \mu
)>m \\ \omega \cdot (\mu '- \nu ')>m
 } } \frac{\bar
\xi^{\mu+\mu'} \xi^{\nu'+\nu}}{ \overline{\xi}_k}\mu_k\bar
c_{\mu\nu\mu'\nu'},
\end{eqnarray}
  where summations are finite and where

\begin{equation}
\label{m.15} c_{\mu\nu\mu'\nu'}:=\langle \Phi_{\mu \nu },R^+_{\mu
'\nu ' }\bar \Phi_{  \nu '\mu '} \rangle \ .
\end{equation}
We further simplify by extracting the main terms. In \eqref{m.14} all
the terms which do not satisfy $ \mu =\nu' =0$ are negligible, see in
particular  the estimate of  \eqref{case1} in Appendix
\ref{prle117}.  In particular, for any of them there is in
\eqref{m.14} a term such that $ \mu =\nu' =0$ which is, clearly,
larger.  In particular all the terms in \eqref{m.14bis} are
negligible (for the proof see the estimate of \eqref{case2} in
  Appendix \ref{prle117}). We ignore all these terms, and proceed in the
discussion. We set $\mathbb{N}_0= \mathbb{N}\cup \{ 0\} $ and we
consider
\begin{equation}
\label{equation:M} M:=\left\{\mu\in \mathbb{N}_0^n\ :\ \mu\cdot
\omega>m \ ,\quad 2 \leq |\mu|\leq 2N+3 \right\} .
\end{equation}
Then, neglecting   all  negligible terms, we write
\begin{equation}
\label{m.13bis} \dot \xi_k =  -\im  \omega_k \xi_k-\im
\frac{\partial Z_0}{\partial
  \bar \xi_k}+\G_{0,k}(\xi)
\end{equation}
where  we set
\begin{equation}
\label{m.13ter} \G_{0,k}(\xi):=\im     \sum_{   \nu  \in M,  \mu \in
M} \frac{\xi^{ \mu }\bar \xi^{ \nu}} { \bar\xi_k} \nu_kc_{0\nu\mu
0}.
\end{equation}
  We focus on \eqref{m.13bis}. Following the idea in
\cite{busper2,sofferweinsten1}, we  apply   normal form theory (in
the form of chapter 5 \cite{Arn83}) in order to further simplify the
system \eqref{m.13bis}.
 We consider a change of variables of the form
\begin{equation}
\label{equation:FGR2}  \eta _j =\xi _j + \Delta_j(\xi)
\end{equation}
which inserted in \eqref{m.13bis} transforms such a system into a
perturbation  (through the small function $\E_j(t)$ defined in
  \eqref{e.es.bis} and estimated in \eqref{e.es.117})   of the system
\begin{equation}
\label{m.17} \dot \eta_k=\Xi_k(\eta,\bar \eta):= -\im \omega_k
\eta_k- \im \frac{\partial Z_0}{\partial
  \bar \eta_k}+\Nc_k(\eta)
\end{equation}
where
\begin{equation}
\label{m.17b}
\Nc_k(\eta):=\im \omega_k\Delta_k(\eta)-\im
\sum_j \left(
\frac{\partial
  \Delta_k}{\partial \eta_j}(\eta)\omega_j\eta_j-  \frac{\partial
  \Delta_k}{\partial \bar\eta_j}(\eta)\omega_j\bar\eta_j
\right)+ \G_{0,k}(\eta).
\end{equation}
The choice
\begin{eqnarray}
\label{m.17bis}
\Delta_j(\xi):= \sum_{\substack{\mu \in M, \, \mu '\in M \\
\omega\cdot(\mu-\nu )\not=0}} \frac{1}{\im \omega \cdot (\mu - \nu) }
\frac{\xi ^{ \mu }\bar {\xi }^ {\nu }}{\bar \xi_j} \nu _j c_{0\nu\mu
0}
\end{eqnarray}
eliminates all non resonant terms from $\N _k$ and reduces it to
\begin{equation}
\label{m.17ter}
\Nc_k(\eta)=\im  \sum_{\substack{\mu \in M, \, \nu\in M \\
\omega\cdot(\mu-\nu ) =0}} \frac{\eta^{\mu}\bar\eta^{\nu}}
{\bar\eta_k} \nu _kc_{0\nu\mu 0}.
\end{equation}
Now we have arrived at the key point of our analysis. Since
$H_{0L}\equiv \sum_{k}\omega_k\left|\eta_k\right|^2$ is a conserved
quantity for the system in which the last term of \eqref{m.17} is
neglected, it is natural to compute the Lie derivative
$\lie_{\Xi}H_{0L}\equiv \sum \omega_j(\bar\eta _j\dot \eta_j+\dot
{\bar \eta }_j\eta _j) $. Notice that we depart here from
\cite{busper2,sofferweinsten1}  and the previous literature, which
rather than at $H_{0L}$, less optimally look   at  $
Q\equiv\sum_{k} \left|\eta_k\right|^2$. The reason for choosing
$H_{0L}$ rather than $Q$ is that $\{ Z _0,H_{0L} \}=0$, while $\{ Z
_0,Q \}=0$  only in the case when all eigenvalues of $-\Delta +V$
are of multiplicity 1. The morale is that with $H_{0L} $  the
multiplicity of the eigenvalues of $-\Delta +V$ is irrelevant in
the argument.  On the other hand, the choice of  $Q$ forces  in the
literature to the hypothesis that the eigenvalues be simple, see
\cite{tsai,zhousigal,cuccagnamizumachi} etc. See also the work in
\cite{zhouweinstein1} in the case of a single multiple eigenvalue
close to the continuous spectrum.

We compute $\lie_{\Xi}H_{0L}$ using Plemelji formula $\frac{1}{x\mp
  i0}=PV\frac{1}{x}\pm\im \pi \delta (x)$, from which one has
$R^{\pm}_{\mu0}=PV(B-\omega\cdot\mu)^{-1}\pm\im \pi \delta
(B-\omega\cdot \mu) $ (where the distributions in $B$ are defined
by means of the distorted Fourier transform associated to $-\Delta
+V$. For the study of positive times, the relevant operator is
  $R^+_{\mu0}$.    Define
\begin{eqnarray} \label{equation:Lambda} \Lambda:=\bigcup_{\mu\in M}\left\{
\omega\cdot \mu \right\} \\ \label{m.32} M_\lambda:=\left\{ \mu\in
M\ :\ \omega\cdot \mu=\lambda\right\} \text{ for $\lambda \in
\Lambda$}
\\
\label{m.32.a} F_\lambda:=\sum_{\mu\in M_\lambda} \bar \eta ^\mu
 \Phi_{0\mu }\ ,\quad B_\lambda:= \pi \delta (B-\lambda).
\end{eqnarray}
Our way to normalize the system leads us to  what follows.

\begin{lemma}
\label{l.dis} The following formula holds:
\begin{equation}
\label{m.33} \lie_{\Xi}H_{0L}=-\sum_{\lambda\in\Lambda}\lambda
\langle F_\lambda;B_\lambda \bar F_{\lambda} \rangle .
\end{equation}
Moreover,  the right hand side is semidefinite negative.
\end{lemma}
\proof We have by \eqref{equation:FGR2} and \eqref{m.15}
\begin{equation} \begin{aligned} &\lie_{\Xi}H_{0L}=-
\Im   \big [ \sum_{\substack{\mu \in M, \, \nu\in M \\
\omega\cdot(\mu-\nu) =0}} \omega \cdot \nu
  \eta^{\mu}\bar\eta^{\nu}
   \langle  \Phi_{ 0\nu},  (B-\omega \cdot \mu -\im  0)
^{-1}\bar \Phi_{  0\mu} \rangle   \big ]   \\&
=-\sum_{\lambda\in\Lambda}\lambda \Im \left [ \langle F_\lambda ,
(B-\lambda -\im 0) ^{-1}  \bar F_\lambda \rangle \right ] .
\end{aligned} \nonumber
\end{equation}
Plemelji formula yields \eqref{m.33}. For $\Psi
_{\lambda}=(B+\lambda )F _{\lambda}$ we have for $k^2=\lambda
^2-m^2$
\begin{eqnarray}
  \langle  F  _{\lambda}, (B-\lambda -\im 0) ^{-1}\bar F_{\lambda}
\rangle   =   \langle    F _{\lambda}, R^{+}_{-\Delta +V}(k^2) \bar
\Psi _{\lambda}\rangle .\nonumber
\end{eqnarray}
The latter is well defined, as stated above in Lemma \ref{lem:m.11}
and proved in  Lemma
\ref{lemma-regularization}  in  Appendix \ref{app.resolv}. By Proposition 2.2 ch. 9 \cite{taylor}
or    Lemma 7 ch. XIII.8 \cite{reedsimon},

\begin{eqnarray}& \Im \left [ \langle    F _{\lambda} , R^{+}_{-\Delta +V}
(k^2)\bar \Psi _{\lambda}\rangle \right ] =\pi \langle    F
_{\lambda} ,
\delta ( -\Delta +V -k^2) \bar \Psi _{\lambda}\rangle = \nonumber \\
& = \frac{k}{16\pi} \int _{|\xi |=k} \widehat{F_{\lambda}} (\xi
)\overline{\widehat{\Psi _{\lambda}} } (\xi ) d\sigma (\xi ) =
\frac{2\lambda k}{16\pi} \int _{|\xi |=k} |\widehat{F_{\lambda}}
(\xi )|^2 d\sigma (\xi ),\nonumber
\end{eqnarray}
where by $\widehat{w}$ we mean the distorted Fourier transform of
$w$ associated to $-\Delta +V$, see Appendix \ref{planewaves}, ch. 9
\cite{taylor} or section XI.6 \cite{reedsimon}.\qed

We will see in subsection \ref{subsec:end} how the structure in
\eqref{m.33}, which continues to hold in the non simplified system,
yields asymptotic stability if we assume the generic conditions
discussed in the next subsection or in (H7). Notice that the sign of
the corresponding term in \cite{cuccagnamizumachi}, see formula
(5.11) \cite{cuccagnamizumachi}, is unclear. Notice that the sign in
(5.11) \cite{cuccagnamizumachi} is nonnegative in the case of 1
eigenvalue, by an indirect argument, see Corollary 4.6
\cite{cuccagnamizumachi}. But here we are interested in the general
case, with many eigenvalues. See also the very complicated argument
in \cite{Gz} to prove the structure \eqref{m.34} in very special
cases (1 eigenvalue with $N=2,3$).

\subsection{The nondegeneracy assumption}
\label{subsection:nondegeneracy}

We are ready to state  the nondegeneracy assumption  mentioned in
the introduction. Specifically, we assume:

\begin{itemize}
\item[(H7)] there exists a positive constant $C$ and a sufficiently
small $\delta _0>0$ such that such that for all $|\eta |<\delta _0$
\begin{equation}
\label{m.34}
\sum_{\lambda\in\Lambda}\lambda \langle F_\lambda;B_\lambda
\bar F_{\lambda} \rangle\geq C \sum_{\mu\in M}|\eta^\mu|^2 \ .
\end{equation}
\end{itemize}
  $M$ and $\Lambda $ are   large sets, so we characterize \eqref{m.34}
  in terms of somewhat smaller sets. Set
\begin{eqnarray} \label{realM} & \widehat{M} =\left\{\mu\in
  M\ :\ \text{$\nu _j \le \mu_j$  $\forall $ $j$ and $\nu \neq \mu$}  \Rightarrow
   \nu \not \in M \right\}
\\ &
\label{realLambda} \widehat{\Lambda}:=\bigcup_{\mu\in
\widehat{M}}\left\{ \omega\cdot \mu \right\}\\&
\widehat{M}_\lambda:=\left\{ \mu\in \widehat{M}\ :\ \omega\cdot
\mu=\lambda\right\} \text{ for $\lambda \in \widehat{\Lambda}$} .
\nonumber
\end{eqnarray}
It is easy to show that (H7) is  equivalent to:
\begin{itemize}
\item[(H7')] For any $\lambda\in\widehat{\Lambda}$ the following matrix is
invertible:
\begin{equation}
\label{m.34a} \left\{\left\langle\bar
\Phi_{\mu0},B_\lambda\Phi_{\mu'0} \right\rangle
\right\}_{\mu,\mu'\in \widehat{M}_\lambda}.
\end{equation}
\end{itemize}

\begin{remark}
\label{stabo} The set $\hat \Lambda$ depends on $m$; $\widehat
M_\lambda$ is piecewise constant in  $m$.
\end{remark}
In the case where $j\neq l$ implies $-\lambda _j^2\neq -\lambda
_l^2$ (this can be easily arranged picking $V(x)$ generic, by
elementary methods in perturbation theory), the assumption (H7) can
be further simplified. Indeed (H5) implies that for any $\lambda\in
\widehat{\Lambda}$ there exists a unique $\mu\in
\widehat{M}_\lambda$. Then (H7') reduces to
\begin{itemize}
\item[(H7'')] For any $\mu\in \widehat{M}$ one has
$\gamma_{\mu}:=\langle \bar \Phi_{\mu0}, B_{\omega\cdot \mu}
\Phi_{\mu0} \rangle\not =0$.
\end{itemize}
We are now ready to give the proof of Proposition \ref{H6}.

\noindent {\bf Proof of Proposition \ref{H6}.} We use equation
\eqref{ze.2} in order to compute the quantities \eqref{m.34a} as
functions of $m$ and of the Taylor coefficients $\beta_l$ of
$\beta$. Set $c=c_\mu =\frac{2^{-\frac{|\mu |}{2}}}{\mu !}$ and
$\Psi_\mu:=B^{-1/2}\varphi^\mu $. Then,
   \eqref{ze.2} implies
\begin{eqnarray}
\nonumber
\gamma_\mu(m,\beta_4,...,\beta_{|\mu|+1}) \\
\label{m.m9}
=
\gamma_\mu(m,\beta_4,...,\beta_{|\mu|},0)+ 2c\beta_{|\mu|+1}
\Re\langle\Phi_{\mu,0}(m,\beta_4,...,\beta_{|\mu|},0),B_{\omega\cdot\mu}
\Psi_\mu \rangle \\
\nonumber +c^2\beta_{|\mu|+1}^2\langle\bar \Psi_\mu ,B_{\omega\cdot
\mu}\Psi_{\mu}\rangle .
\end{eqnarray}
We conclude that either \eqref{m.m9} is independent of
$\beta_{|\mu|+1}$ or there exists at most two values of
$\beta_{|\mu|+1}$ for any choice of $(m,\beta_4,...,\beta_{|\mu| })$
such that \eqref{m.m9} vanishes. We show now that, except for at
most a finite number of values of $m$ in any compact interval,
\eqref{m.m9} depends on $\beta_{|\mu|+1}$. We have, see the proof of
\eqref{m.33},
\begin{equation} \label{m.m10}\langle \bar \Psi_{\mu},
B_{\omega\cdot\mu } \Psi_{\mu} \rangle  =\frac{1}{16\pi }\int _{|\xi
|=\sqrt{ (\omega\cdot\mu ) ^2-m^2}} |  \widehat{\varphi^\mu}  (\xi
)|^2 d\sigma (\xi )\, ,
\end{equation}
  where we are using the distorted Fourier transform
associated to $-\Delta +V$. Since the $\varphi _j (x)$ are smooth
functions decaying like $e^{-|x| |\lambda _j| }$ with all their
derivatives, and $V(x)$ is chosen exponentially decreasing as well,
by Paley Wiener theory applied to the distorted Fourier transform
associated to $-\Delta +V$, the functions $\widehat{ \varphi ^\mu }
(\xi )$ are analytic, see Remark \ref{remark:planewaves}. If the set
where $\widehat{\varphi^\mu} (\xi )=0$ does not contain any sphere,
then the proof is completed. If $\widehat{\varphi^\mu} (\xi )=0$ on
a sphere, say   $|\xi|= a_0$, then, by analyticity,
$\widehat{\varphi^\mu} (\xi )$ does not vanish identically on nearby
spheres. We eliminate values of $m$ such that $\omega(m)\cdot\mu=
a_0 $. Since $\omega(m)\cdot\mu$ is a nontrivial analytic function
this can be obtained by removing at most a finite number of values
of $m$. Repeating the operation for all $\mu\in \widehat M$ (a
finite set) one gets that, apart from a finite set of values of $m$,
the quantity in \eqref{m.m10} is different from 0. Thus removing at
most two values of $\beta_{|\mu|+1}$ for each $\mu\in\widehat M$,
one gets $\gamma_\mu>0$ $\forall \mu\in\widehat M$. \qed

\begin{remark} \eqref{m.m10} with $\mu =3$ and
$\ker (-\Delta +V+\lambda ^2)=\text{span}\{ \varphi \}$ is the
condition necessary in the special case   in \cite{sofferweinsten1}.
If   $\widehat{ \varphi ^3} (\xi )=\widehat{ \varphi ^3 } (|\xi| )$,
then the fact that \eqref{m.m10} is nonzero reduces to $\widehat{
\varphi ^3 } (\sqrt{9\omega ^2-m^2} )\neq 0$, which is the condition
written in (1.8) \cite{sofferweinsten1}.
\end{remark}

\section{Review of linear theory}

We collect here some well known facts needed in the paper. First of
all, for our purposes the following Strichartz estimates for the flat
equation will be sufficient, see \cite{danconafanelli}:

\begin{lemma}\label{lemma-strichartzflat} There is a fixed $C$ such
that for any admissible pair $(p,q) $, see (\ref{admissiblepair}),
we have

\begin{equation}\label{strichartzflat} \|   K'_0(t)
u_0+ K_0(t) v_0 \|
_{L^p_tW^{\frac{1}{q}-\frac{1}{p}+\frac{1}{2},q}_x}\le
 C\| (u_0,v_0) \| _{H^1\times L^2}.
\end{equation}
Furthermore, for any other admissible pair $(a,b)$,
\begin{equation}\label{strichartzflatretarded}
\|  \int _{s<t} K_0(t-s)F(s)ds \|
_{L^p_tW^{\frac{1}{q}-\frac{1}{p}+\frac{1}{2},q}_x}\le
 C \| F \| _{L^{a'}_tW^{\frac{1}{a}-\frac{1}{b} +\frac{1}{2},
 b'}_x},
\end{equation}
where given any $p\in [1,\infty ]$ we set $p'=\frac{p}{p-1}.$
\end{lemma}

We next consider the linearization of (\ref{NLKG}). Notice that
under (H1) for any $k\in \mathbb{N}\cup \{ 0 \}$ and $p\in [1,\infty
]$ the functionals $ \langle
 \cdot ,\varphi _j\rangle  $ are bounded in $W^{k,p}$. Let
$W^{k,p}_c$,  $H^k_c$ if $p=2$, be the intersection of their kernels
in $W^{k,p} $. We recall the following result by \cite{yajima}.
\begin{theorem}\label{yajima} Assume: (H2);
$|\partial _x^\alpha V (x)|\le C \langle x \rangle ^{-\sigma }$ for
$|\alpha |\le k$,  for   fixed $C$ and $\sigma >5$. Consider the
strong limits
\begin{equation}\label{waveoperators} \mathcal{W}_\pm =
\lim _{t\to \pm \infty}e^{\im t (-\Delta +V)}e^{\im t  \Delta  } \,
, \quad \mathcal{Z}_\pm = \lim _{t\to \pm \infty}e^{-\im t \Delta
}e^{\im t ( \Delta -V)}P_c.
\end{equation}
Then $\mathcal{W}_\pm : L^2\to  L^2_c$      are isomorphic
isometries which extend into isomorphisms $\mathcal{W}_\pm :
W^{k,p}\to W^{k,p}_c$ for all $p\in [1,\infty ]$. Their inverses are
$\mathcal{Z}_\pm $. For any Borel function $f(t)$ we have, for a
fixed choice of signs,
\begin{equation}\label{coniugation} f(-\Delta +V)P_c =
\mathcal{W}_\pm f(-\Delta  ) \mathcal{Z}_\pm \, , \quad f(-\Delta
)P_c = \mathcal{Z}_\pm f(-\Delta +V)P_c \mathcal{W}_\pm .
\end{equation}

\end{theorem}

Because of $ \frac{1}{q}-\frac{1}{p}+\frac{1}{2} = \frac{5}{2}
 (\frac{1}{2}-\frac{1}{q})\in [0,5/6 ]$   for all admissible pairs
$(p,q)$, by Theorem \ref{yajima} for $k\le 2$ we have the following
transposition of Lemma \ref{lemma-strichartzflat} to our non flat
case.

\begin{lemma}\label{lemma-strichartz1}
  Set $ K(t)=
 {\sin (t B)}/{B}.$ Then, if we assume (H1)--(H2)
there is a fixed constant $C_0$ such that for any two admissible
pairs $(p,q)$ and $(a,b)$ we have
\begin{equation}\label{strichartz1}\begin{aligned} &\|   K'(t)
u_0+ K(t) v_0 \|
_{L^p_tW^{\frac{1}{q}-\frac{1}{p}+\frac{1}{2},q}_x}\le
 C_0\| (u_0,v_0) \| _{H^1\times L^2}.
\\ &
\|  \int _{s<t} K(t-s)F(s)ds \|
_{L^p_tW^{\frac{1}{q}-\frac{1}{p}+\frac{1}{2},q}_x}\le
 C_0 \| F \| _{L^{a'}_tW^{\frac{1}{a}-\frac{1}{b}+\frac{1}{2}, b'}_x}.
\end{aligned}\end{equation}
\end{lemma}

By   Theorem \ref{yajima}  for $k\le 2$  we have the following
transposition  of the analogous estimates of the flat case, which in
turn are equivalent to Lemma \ref{lemma-strichartzflat}.
\begin{lemma}\label{lemma-strichartz} If we assume (H1)--(H2) there is a fixed constant $C_0$
  such that for any two admissible
pairs $( p,q)$ and $(a,b)$ we have
\begin{equation} \label{strichartz2}\begin{aligned} & \|   e^{-\im tB}P_c u_0 \|
_{L^p_tW^{\frac{1}{q}-\frac{1}{p} ,q}_x}\le
 C_0\|  u_0   \| _{ H^{1/2}}
\\
&     \|  \int _{s<t}  e^{\im (s-t)B}P_c F(s)ds \|
_{L^p_tW^{\frac{1}{q}-\frac{1}{p} ,q}_x}\le
 C _0\| F \| _{L^{a'}_tW^{\frac{1}{a}-
 \frac{1}{b}+1, b'}_x}
.
\end{aligned}
\end{equation}
\end{lemma}

Sketches of proofs of Lemmas \ref{lemma:weighted}  and
\ref{lemma:christkiselev} are in Appendix A.

\begin{lemma}
\label{lemma:weighted} Assume (H1)--(H2) and
  consider $m<a<b<\infty$. Then for any $\gamma>9/2$
  there is a constant $C=C(\gamma )$ such that  we have
\begin{equation}\label{2.6}
\|  e^{ -\im B t}R_{B}(\mu+\im 0) g\|_{H^{-4,-\gamma}_x} \le C\langle
t \rangle ^{-\frac 32} \| g \|_{L^{2,\gamma} _x} \text{ for  any
$\mu \in [a,b]$
  and $t\ge 0$  .}
\end{equation}
\end{lemma}

\begin{lemma}
\label{lemma:christkiselev} Assume (H1)--(H2).Then for any $s>1$
there is a fixed $C_0=C_0(s,a)$ such that for any admissible pair
$(p,q)$ we have
\begin{eqnarray}\label{eq:christkiselev}
\left\|  \int _{0} ^t e^{\im (t'-t)B}P_cF(t') dt'\right \|
_{L^p_tW^{\frac{1}{q}-\frac{1}{p} ,q}_x} \le C_0    \| B^{\frac
12}P_cF\|_{L_t^aL ^{ 2, s }_x}  \end{eqnarray} where for $p>2$ we
can pick any $a\in [1,2]$ while for $p=2$ we pick $a\in [1,2)$.
\end{lemma}

\section{Nonlinear estimates}
We apply Theorem \ref{main} for $r=2N $  (recall $N=N_1$ where
$N_j\omega _j<m <(N_j+1)\omega _j).$   Then we study the solutions of
the Hamilton equations of $H^{(2N)}$ with initial data corresponding
to orginal ones. In particular $f$ and $\xi$ denote the solutions of
such equations.

We will show:
\begin{theorem}\label{proposition:mainbounds} There exist constants
$C >0$ and  $\varepsilon _0>0$ such that, if the initial data in terms
  of the original variables fulfill $\norma{(u_0,v_0)}_{H^1\times
    L^2}\leq \epsilon$, with
$\epsilon \in (0, \varepsilon _0)$, then we have
\begin{eqnarray}
&   \|  f \| _{L^p_t( \mathbb{R},W^{ {1}/{q}- {1}/{p} ,q}_x)}\le
  C \epsilon \text{ for all admissible pairs $(p,q)$} \label{Strichartzradiation}
\\& \| \xi ^\mu \| _{L^2_t(\mathbb{R})}\le
  C \epsilon \text{ for all multi indexes $\mu$
  with  $\omega \cdot \mu >m $} \label{L^2discrete}\\& \| \xi _j  \|
  _{W ^{1,\infty} _t  (\mathbb{R}  )}\le
  C \epsilon \text{ for all   $j\in \{ 1, \dots , n\}$ }
  \label{L^inftydiscrete} .
\end{eqnarray}
\end{theorem}
 Theorem \ref{proposition:mainbounds}
implies \eqref{Strichartz}. The existence of $(u_\pm , v_\pm )$ is
instead a consequence of Lemma \ref{lemma:asymoptoticflatness}
below.

\begin{remark}
\label{r.es.2} By \eqref{boundedenergynorm}  one has $|\xi
|_{L^\infty_t(\R  )} +
 \|
 f  \| _{L^\infty_t(\R ,H^{\frac{1}{2}}_x)}\lesssim \epsilon$.
  Also \eqref{L^inftydiscrete} is an easy consequence of
\eqref{boundedenergynorm} and \eqref{e.9}, so it will be assumed.
\end{remark}
\begin{remark} By the time reversibility of \eqref{NLKG} it
 is not restrictive to prove  Theorem
\ref{proposition:mainbounds} with $\mathbb{R}$ replaced by
$[0,\infty )$. So in the sequel we will consider $t\ge 0$ only.
\end{remark}

\begin{remark}
\label{r.es.1} We have for any bounded interval $I$
\begin{equation}
\label{es.e.1} f\in L^p _t(I,W^{ {1}/{q}- {1}/{p}
  ,q}_x) \text{ for all admissible pairs $(p,q)$} \ .
\end{equation}
This can be seen as follows. $u\in L^\infty _t (\mathbb{R}, H^1_x)$,
  implies $u^3\in L^\infty _t (\mathbb{R}, L^2_x)$ and $\| \beta '( u)\|
_{L^2_x}\le \| u \| _{L^6_x}^3\lesssim \| u \| _{H^1_x}^3$. By Lemma
\ref{lemma-strichartz1} and \eqref{duhamelflat}, this implies $u\in
 L^p _t(I,W^{ {1}/{q}- {1}/{p}
  ,q}_x) $ over any bounded
interval $I$ for any admissible pair $(p,q)$. Then, the estimate
\eqref{def1} implies that the property persists also after the
normalizing transformation.
\end{remark}
We prove Theorem \ref{proposition:mainbounds} by means of a standard
continuation argument, spelled out for example in formulas
(2.6)--(2.8) \cite{sogge}.  We know that $\| f (0)\| _{H^{1/2}}
+|\xi (0) |\le c_0\epsilon $. We can consider a fixed constant $C_3$
valid simultaneously for Lemmas
\ref{lemma-strichartz}--\ref{lemma:christkiselev}. Suppose that the
following estimates hold

\begin{eqnarray}
&   \|  f \| _{L^p_t([0,T],W^{ {1}/{q}- {1}/{p} ,q}_x)}\le
  C _1\epsilon \text{ for all admissible pairs $(p,q)$} \label{4.4a}
\\& \| \xi ^\mu \| _{L^2_t([0,T])}\le
 C_2 \epsilon \text{ for all multi indexes $\mu$
  with  $\omega \cdot \mu >m $} \label{4.4}
\end{eqnarray}
 for fixed large multiples $C_1$, $C_2$ of $c_0C_3$. Then we will
 prove that, for $\epsilon $ sufficiently small independent of $T$,
 \eqref{4.4a} and \eqref{4.4} imply the same estimate but with $C_1$,
 $C_2$ replaced by $C_1/2$, $C_2/2$.  Then \eqref{4.4a} and
 \eqref{4.4} hold with $[0,T]$ replaced by $[0,\infty )$.

\subsection{Estimate of the continuous variable $f$}
 Consider $H^{(2N)}=H_L+Z^{(2N)}+ \resto ^{(2N)}$. We set $Z=Z^{(2N)}$ and
 $\resto= \resto ^{(2N)}$. Then we have

\begin{eqnarray} \label{equation:f}
& \im \dot f-  Bf=  \nabla _{\bar f} Z_1+  \nabla _{\bar f}
\resto \end{eqnarray}

\begin{lemma}\label{lemma:bound remainder}
Assume \eqref{4.4a}, and \eqref{4.4}, and fix  a large $s>0$.  Then
there is a   constant $C=C(C_1,C_2)$ independent of $\epsilon$ such
that the following is true: we have $ \nabla _{\bar f} \resto =R_1+
R_2$ with
\begin{equation}
\label{bound1:z1} \| R_1 \| _{L^1_t([0,T],H^{\frac{1}{2} }_x)}+\|
B^{\frac 12}P_c R_2 \| _{L^{2
\frac{N+1}{N+2}}_t([0,T],L^{2,s}_x)}\le C(C_1,C_2) \epsilon^2.
\end{equation}
\end{lemma}
\proof For $d\le 1 $ and arbitrary fixed $s$  we have $\nabla
_{\bar f}\resto _d \in H^{\frac{1}{2},s} $.
By   (iii0--iii1) and Theorem \ref{main}
$$\| \nabla _{\bar f}\resto _0  \|_{ H^{\frac{1}{2},s} } + \|
\nabla _{\bar f}\resto _1  \|_{ H^{\frac{1}{2},s} }\le C|\xi  |
^{2N+3}.$$ Hence by \eqref{4.4}  and Remark \ref{r.es.2}

\begin{equation}
\label{es.e.3} \norma{\nabla _{\bar f}(\resto _0 +\resto
_1)}_{L^1_t([0,T], H^{\frac{1}{2} }_x)} \lesssim  \| |\xi |^{N+1} \|
^{2 }_{L^{2 }_t} \| \xi \|  _{L^{\infty}_t} \le C_2^2C \epsilon ^{
3}.
\end{equation}
  $\nabla_{\bar
  f}\resto _d$ with $d\le 1 $ is  absorbed in $R_1$. For $d=2,3$ we have
  \begin{equation}
\label{schematic11}\begin{aligned} &\nabla_{\bar
  f}\resto _d=\frac{d}{\sqrt{2}}\,  B^{-\frac{1}{2}} ( F_d(x,z,B^{-\frac{1}{2}} f(t,\cdot )) U^{d-1}(t,\cdot ))+ \\& +\frac{1}{\sqrt{2}} B^{-\frac{1}{2}} ( \partial _w F_d(x,z,B^{-\frac{1}{2}} f(t,\cdot )) U^{d }(t,\cdot ))+\\& + \nabla _{\overline{g}}\left (  \int_{\R^3}  F_d  (x,\xi , g ,B^{-\frac{1}{2}}f(t,x))[U(t,x)]^ddx\right ) _{g=f} .
\end{aligned}
\end{equation}
Similarly for $(\xi ',f')=\Tr (\xi , f)$ and $U'= \frac{1}{\sqrt{2B}} (f'+\overline{f'})$ we have

\begin{equation}
\label{schematic12}\begin{aligned} &\nabla_{\bar
  f}\resto _4=2\sqrt{2}   B^{-\frac{1}{2}} ( F_4(x,\xi ', U'(t,\cdot )) U^{3}(t,\cdot ))+ \\& +\frac{1}{\sqrt{2}}B^{-\frac{1}{2}} ( \partial _Y F_4(x,\xi ', U'(t,\cdot )) U^{4 }(t,\cdot ))+\\& + \sum _{j=1}^{n} \int_{\R^3}\partial  _{\xi '_j}  F_4  (x, \xi ', U'(t,x) ) [U(t,x)]^4dx\, \,  \nabla   _{ \overline{f}} \xi _j' \\& + \sum _{\mu \nu} \int_{\R^3}\partial  _{Y}  F_4  (x,  \xi ', U'(t,x)  )  \Psi  _{\mu \nu}  (x)  [U(t,x)]^4dx\, \, \nabla_{\bar
  f} G  _{\mu \nu}  (z),
\end{aligned}
\end{equation}
$G _{\mu,\nu}$ as in Lemma \ref{lie_trans},
$\Psi _{\mu\nu}(x)\in\Sc(\R^3,\C)$ and   $Y$ as in (3) Lemma \ref{lemma:H_P}.

The sums of the  contributions from the first two lines of   \eqref{schematic11}--\eqref{schematic12} are   schematically of the form
\begin{equation}
\label{schematic1}\begin{aligned} &
         B^{-\frac{1}{2}} \left [ \left ( \Phi _1 (x,z)  B^{-\frac{1}{2}}f  \right )
    +   \left ( \Phi _2 (x) (B^{-\frac{1}{2}}f )^2\right )
 + f^3\right ] ,
\end{aligned}
\end{equation} with a   $\Phi _2 \in H^{k,s}(\R ^3, \C) $  and with
$\Phi _1 (x,z)\in C^\infty ( \U ^{-k,-s}_z, H^{k,s}(\R _x^3, \C) )$
such that $\| \Phi _1 (x,z) \| _{H^{k,s}}\le C\norma{z}_{\Ph
^{-k,-s} }.$   $R_2$ is formed
by the first term in \eqref{schematic1}, while all the rest can be
absorbed in $R_1$. The last line of  \eqref{schematic11} and the last two lines of  \eqref{schematic12}
are  absorbed in $R_1$.
Let us start with the terms forming $R_1$.

By Theorem \ref{yajima}, using the wave operator $\mathcal{Z}_+$ in
\eqref{waveoperators}, we have
\begin{equation} \label{estimate1}\begin{aligned} &
  \|   B^{-\frac{1}{2}}\left ( \Phi _2(x) (B^{-\frac{1}{2}}f )^2
  \right )
   \| _{L^1_tH^{\frac 12}_x}\lesssim \|  \mathcal{Z}_+ B^{-\frac{1}{2}}\left ( \Phi _2(x) (B^{-\frac{1}{2}}f )^2
  \right )
   \| _{L^1_tH^{\frac 12}_x}  \\& =\|  (-\Delta +m^2)^{-\frac{1}{2}}
    \mathcal{Z}_+\left ( \Phi _2(x) (B^{-\frac{1}{2}}f )^2
  \right )
   \| _{L^1_tH^{\frac 12}_x}
   \\&  \lesssim  \|    \left ( \Phi _2(x) (B^{-\frac{1}{2}}f )^2
  \right )
   \| _{L^1_tL^{ 2}_x} \lesssim
     \|   B^{-\frac{1}{2}}f  \| _{L^2_tL^{6}_x}^2 \lesssim
     \|  \mathcal{Z}_+ B^{-\frac{1}{2}}f  \| _{L^2_tL^{6}_x}^2\\& =
     \|  (-\Delta +m^2)^{-\frac{1}{2}}\mathcal{Z}_+ f  \| _{L^2_tL^{6}_x}^2
     \approx \|
      \mathcal{Z}_+ f  \| _{L^2_tW^{-1/2,6}_x}^2\lesssim \|
       f  \| _{L^2_tW^{-1/2,6}_x}^2 \\&
     \lesssim
     \|   f  \| _{L^2_tW^{-1/3, 6}_x}^2 \le C_1^2 \epsilon ^2
,
       \end{aligned}
\end{equation}
where in the last line we used \eqref{4.4}.  Proceeding similarly,
    by Remark \ref{r.es.2} and (H6),

\begin{equation} \label{estimate2}\begin{aligned} &
\|  B^{-\frac{1}{2} }     (B^{-\frac{1}{2}}f) ^3  \| _{L^1_tH^{\frac
12}_x}\lesssim \|   (B^{-\frac{1}{2}}f) ^3 \|
_{L^1_tL^{2}_x}\lesssim
   \|   B^{-\frac{1}{2}}f  \| _{L^\infty_tL^{6}_x}
     \|   B^{-\frac{1}{2}}f  \| _{L^2_tL^{6}_x}^2  \\&
\lesssim  \|
 f  \| _{L^\infty_tW^{-\frac{1}{2},6}_x}
     \|   f  \| _{L^2_tW^{-1/3, 6}_x}^2 \lesssim  \|
 f  \| _{L^\infty_tH^{\frac{1}{2}}_x}
     \|   f  \| _{L^2_tW^{-1/3, 6}_x}^2 \lesssim
      C_1^2 \epsilon ^3.
       \end{aligned}
\end{equation}
Looking at   the third line of \eqref{schematic11} we have
\begin{equation} \label{estimate23}\begin{aligned} & \| \nabla _{\overline{g}} \int_{\R^3} F_d  (x,\xi , g ,B^{-\frac{1}{2}}f(t,x))_{g=f}[U(t,x)]^ddx \|  _{L^1_tH^{ 1/2}_x}=\\& \| \sup _{\| \psi  \|  _{ H^{ -\frac 1 2}_x} =1} \int_{\R^3}   d _{\overline{g}} F_d  (x,\xi , g ,B^{-\frac{1}{2}}f(t,x))_{g=f}[\psi ] [U(t,x)]^ddx \| _{L^1_t}.
 \end{aligned}
\end{equation}
For $d=2$ by \eqref{estimate1} and  by  \eqref{eq:coeff F1}  the rhs of  \eqref{estimate23} is

\begin{equation} \label{estimate24}\begin{aligned} & \le C \sup _{\| \psi  \|  _{ H^{ -\frac 1 2}_x} =1} \| d _{\overline{g}} F_2  (x,\xi , g ,B^{-\frac{1}{2}}f(t,x))_{g=f}[\psi ] \| _{L^{\frac{3}{2} } _x}  \| B^{-\frac{1}{2}}f \| _{L^2_tL^6_x}^2\le  C C_1^2\epsilon ^2. \end{aligned}
\end{equation}
For  $d=3$ by \eqref{estimate1}   the rhs of  \eqref{estimate23} is similarly $\le $\begin{equation} \label{estimate241}\begin{aligned} &   C \sup _{\| \psi  \|  _{ H^{ -\frac 1 2}_x} =1} \| d _{\overline{g}} F_3  (x,\xi , g ,B^{-\frac{1}{2}}f(t,x))_{g=f}[\psi ] \| _{L^{2 } _x}  \| (B^{-\frac{1}{2}}f)^3 \| _{L^1_tL^2_x} \le  C C_1^3\epsilon ^3. \end{aligned}
\end{equation}
 We have by \eqref{estimate2}
 \begin{equation} \label{estimate25}\begin{aligned} &     \|
\int_{\R^3} \partial  _{\xi '_j}  F_4  (x, \xi ', U'(t,x) ) [U(t,x)]^4 dx\| _{L^1_t }  \,  \| \nabla   _{ \overline{f}} \xi _j' \| _{H^{ \frac{1}{2} } _x}\\& \le  C \|  B^{-\frac{1}{2}}f  \| _{L^\infty _tL^2_x}  \| (B^{-\frac{1}{2}}f)^3 \| _{L^1_tL^2_x} \le  C C_1^4\epsilon ^4  \end{aligned}
\end{equation}
and
\begin{equation} \label{estimate26}\begin{aligned} &     \|
\int_{\R^3} \partial  _{Y}  F_4  (x,  \xi ', U'(t,x)  ) \Psi  _{\mu \nu}  (x) [U(t,x)]^4dx \| _{L^1_t }  \,  \| \nabla_{\bar
  f} G  _{\mu \nu}  (z) \| _{H^{ \frac{1}{2} } _x}\\& \le  C \|  B^{-\frac{1}{2}}f  \| _{L^\infty _tL^2_x}  \| (B^{-\frac{1}{2}}f)^3 \| _{L^1_tL^2_x} \le  C C_1^4\epsilon ^4 . \end{aligned}
\end{equation}

Collecting in $R_1$ all terms estimated in
\eqref{es.e.3} and   \eqref{estimate1}--\eqref{estimate25} yields the   estimate for
$R_1$. Let $R_2$ be a sum of terms of the form $ \xi
B^{-\frac{1}{2}}\left ( \Phi _1 (x )  B^{-\frac{1}{2}}f\right )$.
Then, proceeding as for \eqref{estimate1}--\eqref{estimate2} and by
\eqref{4.4a} and \eqref{4.4}

\begin{equation} \label{estimate0} \begin{aligned} &
\|   \xi
   B^{-\frac{1}{2}}\left ( \Phi _1   B^{-\frac{1}{2}}f  \right )
   \| _{L^{2\frac{N+1}{N+2}}_tH^{\frac 12,s}_x}  \lesssim  \|  \xi
  P_c\left ( \Phi _1   B^{-\frac{1}{2}}f  \right )
   \| _{L^{2\frac{N+1}{N+2}}_tL^{2,s}_x}
\\&
\lesssim
    \| \xi \| _{L^{2N+2}  _t} \|   B^{-\frac{1}{2}}f
    \| _{L^2_tL^{6}_x}
     \le C_1 \epsilon
     \|   f  \| _{L^2_tW^{-1/3, 6}_x}  \le C_2 C_1   \epsilon ^2.
        \end{aligned}
\end{equation}

  \qed
\begin{remark}
\label{rem.es.3} By \begin{eqnarray} \label{pri2scipio} &
   | \nabla _{\bar \xi }\resto  |\lesssim   |\xi |^{2N+3}
   +  |\xi |^{2N+2} \|  B^{-\frac{1}{2}} f \| _{L^{2,-s}_x}\\ &+
   \|  B^{-\frac{1}{2}} f \|^2 _{L^{2,-s}_x}+
   \|   B^{-\frac{1}{2}} f  \| _{L^{2,-s}_x}^{\frac{3}{2}}
   \|   B^{-\frac{1}{2}} f  \| _{L^{6}_x}^{\frac{3}{2}};
  \nonumber
\end{eqnarray}
and by  the same method as above one can   prove for a fixed
$C$
\begin{equation}
\label{e.es.6} \norma{\partial_{\bar\xi}\resto }_{L^1_t}\leq C
C_1(C_2+C_1+C_1^2) \epsilon^2.
\end{equation}
  One also has the
  easier estimate for   fixed
$C$ and $C_0$
\begin{equation}\label{bound1:z2} \|
\int _0^t  e^{ \im B(s-t)}\nabla _{\bar f}Z_1  \|
_{L^p_tW^{\frac{1}{q}-\frac{1}{p} ,q}_x}\le C_0 \| \nabla _{\bar f}
Z_1  \| _{L^2_tW^{\frac{1}{3}+\frac{1}{2}  , \frac{6}{5}}_x}\le C
C_0C_2\epsilon .
\end{equation}
The important fact is that \eqref{bound1:z2} is independent of
$C_1$.
\end{remark}

\begin{proposition}\label{Lemma:conditional4.2} Assume \eqref{4.4a}
  and \eqref{4.4}. Then there exist constants $K_1 $    and   $C=C(C_1,C_2)$
   such that, if
  $C(C_1,C_2) \epsilon <C_0$, with $C_0$ the constant in Lemma
  \ref{lemma-strichartz}, then we have
\begin{eqnarray}
&   \|  f \| _{L^p_t([0,T],W^{ {1}/{q}- {1}/{p} ,q}_x)}\le
  K_1  \epsilon \text{ for all admissible pairs $(p,q)$}\ .
  \label{4.5}
\end{eqnarray}
\end{proposition}
 \proof Using Lemma \ref{lemma:bound remainder}
 we write \begin{eqnarray} \label{duhamel:f}
f= e^{-\im Bt}f(0)
       -\im \int _0^t  e^{ \im B(s-t)}\nabla
_{\bar f}Z  ds
       - \im \sum _{j=1}^{2}\int _0^t  e^{ \im B(s-t)}
   P_c R_j ds .
\end{eqnarray}
 By \eqref{strichartz2} for $(a,b)=(\infty ,2)$  and \eqref{bound1:z1}
\begin{equation}  \label{eq:cond 1} \begin{aligned} &
 \|  \int _0^t  e^{ \im B(s-t)}
    R_1 ds  \| _{L^p_t([0,T],W^{ \frac {1} {q}- \frac{1} {p} ,q}_x)}
    \le C \| R_1 \| _{L^1_t([0,T],H^{\frac{1}{2} }_x)}
    \le C(C_1,C_2) \epsilon^2.
\end{aligned}
\end{equation}
Similarly, by \eqref{eq:christkiselev} and \eqref{bound1:z1}, we get
for $s>1$
\begin{equation}  \label{eq:cond 2} \begin{aligned} &
 \|  \int _0^t  e^{ \im B(s-t)}
   P_c R_2 ds  \| _{L^p_t([0,T],W^{ \frac {1} {q}- \frac{1} {p}  ,q}_x)}
      \le   {C} \|  \sqrt{B} P_cR_2 \| _{L^{2\frac{N+1}{N+2}}_t([0,T],
  L^{2, s }_x)} \\&  \le C(C_1,C_2) \epsilon^2.
\end{aligned}
\end{equation}
Then the proof is obtained by \eqref{eq:cond 1}--\eqref{eq:cond 2},
by \eqref{bound1:z2} and by \begin{equation}\| e^{-\im Bt}f(0) \|
_{L^p_t(\R,W^{ \frac {1} {q}- \frac{1} {p}  ,q}_x)} \le C_0 \| f(0)
\| _{H^{\frac{1}{2}}} \le K_0\epsilon, \nonumber \end{equation}
which follows by \eqref{strichartz2}. \qed

We end this subsection by proving asymptotic flatness of $f$ if Theorem
\ref{proposition:mainbounds} holds.

 \begin{lemma}
\label{lemma:asymoptoticflatness} Assume
Theorem
\ref{proposition:mainbounds}.
Then there exists $f_+\in H^{\frac{1}{2}}_x$ such that
\begin{eqnarray}\label{eq:asymoptoticflatness} \lim _{t\to \pm
\infty}\left\|  f(t) -e^{-\im Bt}f_+ \right\| _{H^{\frac{1}{2}}_x}
=0.
\end{eqnarray}
\end{lemma}
\proof We have
 $$  e^{\im tB} f(t)=f(0) -\im \int _0^t e^{\im sB}
  \nabla _{\bar f}(Z_1+ \resto )ds $$   and so for
 $t_1<t_2$
 $$  e^{\im t_2B} f(t_2)-
  e^{\im t_1B} f(t_1)=-\im \int _{t_1}^{t_2} e^{\im t'B}
 \nabla _{\bar f}(Z_1+ \resto )  dt'. $$ By
 Lemmas  \ref{lemma-strichartz}, \ref{lemma:christkiselev} and \ref{lemma:bound
 remainder} and by \eqref{bound1:z2},
we get for $t_1\to \infty$ and $t_1<t_2$ \begin{equation}
\begin{aligned} & \| e^{\im t_2B} f(t_2)-
  e^{\im t_1B} f(t_1) \|
_{H^{\frac{1}{2}}_x}
 = \| \int _{t_1}^{t_2} e^{\im t'B}
 \nabla _{\bar f}(Z_1+ \resto )
 dt' \| _{H^{\frac{1}{2}}_x}
\leq \\& \| R_1\| _{ L^1_t([t_1,t_2],H^{\frac{1}{2}}_x) } + \|
\sqrt{B}P_c R_2\| _{L^{2\frac{N+1}{N+2}}_t([t_1,t_2],L^{{2},s}_x) }+ \|
 \nabla _{\bar f}Z_1\| _{L^{2}_t([t_1,t_2],W^{\frac{5}{6},\frac{6}{5} }_x)  }
 \to 0  .\end{aligned} \nonumber
\end{equation}
Then   $f_+=\lim _{t\to \infty}e^{\im tB} f(t) $ satisfies the
desired properties. \qed

 Lemma \ref{lemma:asymoptoticflatness} implies the existence of the
 $(u_+ , v_+ )$ and their properties in Theorem \ref{theorem-1.1}.
\subsection{Estimate of   $g$}

Consider  the $g$ defined in  \eqref{m.6}, \eqref{m.8},
\eqref{m.11}. If $f$, $\xi$ satisfy the Hamilton equations of
\eqref{eq:bir1}, then $g$ satisfies
\begin{equation} \label{e.es.101}\begin{aligned} &
 \im \dot g-Bg=  \nabla_{\bar f}\resto +
  \sum _k \left [  \partial _{ \xi_k} \bar Y
\partial _{  \bar \xi_k}\left( Z+\resto \right )
-\partial _{  \bar \xi_k} \bar Y \partial _{ \xi_k}\left( Z+\resto
\right ) \right] .
\end{aligned}
\end{equation}
We have:

\begin{lemma}\label{lemma:bound g} Assume \eqref{4.4a} and
\eqref{4.4}. Fix $s>9/2$. Then, there are constants $\epsilon _0>0$
and $C>0$ such that, for $\epsilon \in (0,\epsilon _0) $ and   for
$C_0$ the constant in Lemma \ref{lemma-strichartz}, we have
\begin{equation} \label{bound:auxiliary}\| g
\| _{L^2_t([0,T],H^{-4,-s}_x)}\le C_0 \epsilon + C\epsilon ^2
.\end{equation}
\end{lemma}
\proof We can apply Duhamel formula and write
\begin{eqnarray}\label{duhamel:g} g(t)=e^{-\im Bt}g(0)
-\im  \int _0^te^{\im B(t'-t)} [\nabla_{\bar f}\resto  +
\text{second term rhs\eqref{e.es.101}}] dt'.
\end{eqnarray}

First of all we prove $\| e^{-\im Bt}g(0) \| _{L^2_tH^{-4,-s}_x}\le
C_0 \epsilon + O(\epsilon ^2)$. To this end recall that $g(0) =
f(0)+\bar Y(0)$. By Schwarz and  Strichartz inequalities (see Lemma
\ref{lemma-strichartz}) we have
\begin{equation} \| e^{-\im Bt}f(0) \| _{L^2_tH^{-4,-s}_x}\lesssim
\| e^{-\im Bt}f(0) \| _{L^2_tW^{- \frac{1}{3},6}_x}\le C_0\epsilon .
\nonumber\end{equation} The estimate of $\| e^{-\im Bt}\bar Y(0) \|
_{L^2_tH^{-4,-s}_x} $ follows from
\begin{eqnarray} \| e^{-\im Bt}\xi ^ \mu (0) \bar {\xi} ^\nu (0)
 R_{\mu \nu}^+\overline{ \Phi } _{\nu \mu } \| _{L^2_tH^{-4,-s}_x}
 \lesssim |\xi ^ \mu (0) \overline{\xi} ^\nu (0) | \| \overline{ \Phi
 } _{\nu \mu } \| _{ L^{2, s}_x}\lesssim \epsilon ^{|\mu +\nu |}
 ,\nonumber\end{eqnarray} which in turn follows from Lemma
 \ref{lemma:weighted}. We have by Lemma \ref{lemma:bound remainder} and by the proof of Lemma \ref{Lemma:conditional4.2},

 \begin{equation}  \begin{aligned} & \|  \int _0^te^{\im B(t'-t)}  \nabla_{\bar f}\resto  \| _{L^2_tH^{-4,-s}_x} \le \|  \int _0^te^{\im B(t'-t)}  \nabla_{\bar f}\resto  \| _{L^2_tW^{- \frac{1}{3},6}_x}\le C(C_1,C_2)\epsilon ^2.
 \end{aligned}\nonumber
\end{equation}
  The second term in the rhs of  \eqref{e.es.101} contributes  through  various terms
 to \eqref{duhamel:g}. We consider the main ones (for the others the
 argument is simpler). Consider in particular contributions from
 $Z_0$. For $\mu _j\neq 0$ we have by Lemma \ref{lemma:weighted}

\begin{eqnarray} &
\| \int _0^te^{\im (t'-t)B} \frac{{\xi} ^{ {\mu}  }  \bar \xi ^ {
\nu }}{\xi _j}
      \partial _{ \overline{{\xi}} _j}  Z _0
   R_{\mu \nu}^+
 \overline{ \Phi } _{\nu \mu }  dt' \| _{L^2_tH^{-4,-s}_x}\le C
   \|  \frac{{\xi} ^{ {\mu}  }  \bar \xi ^ {
\nu }}{\xi _j}
      \partial _{ \overline{{\xi}} _j}  Z _0 \| _{L^2 _t}
     \| \overline{ \Phi } _{\nu \mu }    \| _{ L^{2, s}_x}   .
    \nonumber\end{eqnarray}
 We need to show

\begin{eqnarray} \label{countingfactors1} &
\| \frac{{\xi} ^{ {\mu}  }  \bar \xi ^ { \nu }}{\xi _j}
      \partial _{ \overline{{\xi}} _j}  Z _0  \| _{L^2 _t}
      =O(\epsilon ^2).
  \end{eqnarray}
By \eqref{m.3} and \eqref{m.11} we have
\begin{eqnarray} \label{countingfactors2}
\omega \cdot (  {\mu} -\nu  ) >m  .\end{eqnarray} Let $\xi ^\alpha
\bar {\xi}^\beta $ be a generic monomial of $Z_0$. The nontrivial
case  is $\beta _j\neq 0$. Then $\partial _{\overline{\xi} _j}(\xi
^\alpha \bar {\xi}^\beta ) =\beta _j\frac{\xi ^\alpha \bar {\xi}^{
{\beta} }}{\overline{\xi}_j}.$ By Definition \ref{d.1} we have
$\omega \cdot (\alpha -\beta ) =0 $, and by Remark \ref{alfabeta},
$|\alpha|=|\beta|\geq 2$. Thus in particular one has
\begin{equation}
\label{ciao}
\omega\cdot \alpha\geq \omega_j\Longrightarrow
\omega\cdot(\mu+\alpha)-\omega_j >m\ .
\end{equation}
So, by remark 7.2 and (7.6), the following holds
\begin{eqnarray}\label{countingfactors4}
\| \frac{{\xi} ^{ {\mu}  }  \bar \xi ^ { \nu }}{\xi _j}
       \frac{\xi ^\alpha \bar {\xi}^{
{\beta} }}{\overline{\xi}_j}  \| _{L^2 _t} \le \| \frac{  \xi ^ {
\nu }  {\xi}^{ {\beta} }}{{\xi}_j}  \| _{L^\infty  _t}\|
\frac{{\xi} ^{  \mu   }
        \xi ^\alpha  }{ {\xi}_j}  \| _{L^2 _t}\leq
C_2C\epsilon^{|\nu|+|\beta|}\leq CC_2\epsilon^2,
  \end{eqnarray}
  where we used $|\xi _l|=|\bar \xi _l|$.
  This completes the proof of Lemma \ref{lemma:bound g}. \qed

\subsection{Estimate of the discrete variables $\xi $}
\label{discrete}

We now return to discrete variables.

\begin{lemma}
\label{eqforeta}
Let $(\xi(t),f(t))$ be a solution of the Hamilton equations of $H^{(2N)}$
and let $(\eta(t),g(t))$ be the corresponding solution defined throgh
\eqref{equation:FGR2} and \eqref{m.6}, then one has
\begin{equation}
\label{e.es.112} \dot\eta_j=-\im \omega_j\eta_j-\im \frac{\partial
    Z_0}{\partial\bar\xi_j}(\eta)+ \N_j (\eta)+\E_j(t)
\end{equation}
where $\N _j$ is defined by \eqref{m.17ter}, and the remainder
$\E_j$ is given by
\begin{eqnarray}
\label{e.es.bis}
 \E_j(t) &:=& \G_{1,j}(\xi)-\im \left\langle \frac{\partial
  G}{\partial \bar \xi_x}(\xi); g \right\rangle
-\im \left\langle \frac{\partial
 \bar  G}{\partial \bar \xi_j}(\xi); \bar g \right\rangle-
\im \frac{\partial \resto^{(2N)}}{\partial
\bar \xi_j}(\xi,f)
\\
\nonumber
&-&\im \sum_{k}\left[\frac{\partial
 \Delta_j}{\partial  \xi_k} \left(  \frac{\partial
 Z^{(2N)}}{\partial \bar \xi_k}+\frac{\partial
 \resto^{(2N)}}{\partial \bar \xi_k}    \right)-
\frac{\partial
 \Delta_j}{\partial \bar  \xi_k} \left(  \frac{\partial
 Z^{(2N)}}{\partial \xi_k}+\frac{\partial
 \resto^{(2N)}}{\partial  \xi_k}    \right)\right]
\\
\nonumber
&+& \left(\N_j(\xi)-\N_j(\eta)
-\im \frac{ \partial Z_0}{\partial \bar
  \xi_j} (\xi)+\im \frac{ \partial Z_0}{\partial \bar
  \xi_j} (\eta) \right)\ ,
\end{eqnarray}
and
\begin{equation}
\label{nu.1} \G_{1,k}(\xi):=
\text{\eqref{m.14}}+\text{\eqref{m.14bis}} -\G_{0,k}(\xi)\ .
\end{equation}
\end{lemma}
\proof First we write the equation for $\xi$. It is convenient to
have in mind the expression in terms of $(\xi, f)$ and an expression
involving also the $g$ variables, namely
\begin{eqnarray}
\label{perxi.1}
\dot \xi_j=-\im \omega_j \xi_j - \im \frac{\partial Z^{(2N)}}{\partial
\bar \xi_j}(\xi,f)-\im \frac{\partial \resto^{(2N)}}{\partial
\bar \xi_j}(\xi,f)
\\
\nonumber
= -\im \omega_j \xi_j - \im \frac{\partial Z_0}{\partial
\bar \xi_j}(\xi,f)- \G_{0,j}(\xi)+L^{(1)}_j(\xi,f,g)\ ,
\end{eqnarray}
where we defined
\begin{equation}
\label{perxi.2}
L^{(1)}_j(\xi,f,g):= \G_{1,j}(\xi)-\im \left\langle \frac{\partial
  G}{\partial \bar \xi_x}(\xi); g \right\rangle
-\im \left\langle \frac{\partial
 \bar  G}{\partial \bar \xi_j}(\xi); \bar g \right\rangle-
\im \frac{\partial \resto^{(2N)}}{\partial
\bar \xi_j}(\xi,f)\ .
\end{equation}
Here and in the rest of the proof, the terms denoted by capital l will
be included in the remainder.

Introducing the variables $\eta$, we have
\begin{equation}
 \begin{aligned} &
\dot \eta_j= \sum_{k} \left(\delta_{jk}+ \frac{\partial
 \Delta_j}{\partial \xi_k} \right)\dot \xi_k+\frac{\partial
 \Delta_j}{\partial \bar \xi_k}\dot\bar\xi_k
\\ &
=\dot \xi_j+ \frac{\partial
 \Delta_j}{\partial \xi_k}\dot \xi_k + \frac{\partial
 \Delta_j}{\partial \bar \xi_k}\dot\bar\xi_k
\\ &
= \dot \xi_j -\im\sum_{k} \omega_k\left(\xi_k\frac{\partial
 \Delta_j}{\partial  \xi_k}-\bar\xi_k\frac{\partial
 \Delta_j}{\partial \bar \xi_k}   \right) +L^{(2)}_j(\xi,f)
 \end{aligned}\nonumber
\end{equation}
where $L^{(2)}_j(\xi,f):=$
\begin{equation}
  -\im \sum_{k}\left[\frac{\partial
 \Delta_j}{\partial  \xi_k} \left(  \frac{\partial
 Z^{(2N)}}{\partial \bar \xi_k}+\frac{\partial
 \resto^{(2N)}}{\partial \bar \xi_k}    \right)-
\frac{\partial
 \Delta_j}{\partial \bar  \xi_k} \left(  \frac{\partial
 Z^{(2N)}}{\partial \xi_k}+\frac{\partial
 \resto^{(2N)}}{\partial  \xi_k}    \right)\right]\ . \nonumber
\end{equation}
Then using the other form of the equations for $\xi$, we have
\begin{eqnarray*}
\dot \eta_j= -\im \omega_j\xi_j-\im \frac{ \partial Z_0}{\partial \bar
  \xi_j} (\xi)+\G_{0,j}(\xi)
\\
-\im\sum_{k} \omega_k\left(\xi_k\frac{\partial
 \Delta_j}{\partial  \xi_k}-\bar\xi_k\frac{\partial
 \Delta_j}{\partial \bar \xi_k}   \right)
\\
+L^{(1)}_j(\xi,f,g)+L^{(2)}_j(\xi,f)\ .
\end{eqnarray*}
Insert now in the first term at r.h.s $\xi_j=\eta_j-\Delta_j(\xi)$.
Thus we get
\begin{eqnarray*}
\dot \eta_j= -\im \omega_j\eta_j+\im \omega_j\Delta_j(\xi)
-\im \frac{ \partial Z_0}{\partial \bar
  \xi_j} (\xi)+\G_{0,j}(\xi)
\\
-\im\sum_{k} \omega_k\left(\xi_k\frac{\partial
 \Delta_j}{\partial  \xi_k}-\bar\xi_k\frac{\partial
 \Delta_j}{\partial \bar \xi_k}   \right)
\\
+L^{(1)}_j(\xi,f,g)+L^{(2)}_j(\xi,f)\ ,
\end{eqnarray*}
which, recalling the definition \eqref{m.17b} of $\N_j$, takes the form
\begin{eqnarray*}
\dot \eta_j= -\im \omega_j\eta_j+\N_j(\xi)
-\im \frac{ \partial Z_0}{\partial \bar
  \xi_j} (\xi)
\\
+L^{(1)}_j(\xi,f,g)+L^{(2)}_j(\xi,f)
\\
= -\im \omega_j\eta_j+\N_j(\eta)
-\im \frac{ \partial Z_0}{\partial \bar
  \xi_j} (\eta)
\\
+ \left(\N_j(\xi)-\N_j(\eta)
-\im \frac{ \partial Z_0}{\partial \bar
  \xi_j} (\xi)+\im \frac{ \partial Z_0}{\partial \bar
  \xi_j} (\eta) \right)+L^{(1)}_j(\xi,f,g)+L^{(2)}_j(\xi,f)\ .
\end{eqnarray*}
Defining  $\E_j$ as the last line of this formula one has the result. \qed

We have:
\begin{lemma}
\label{l.es.11} There is a fixed $C$ such that for
  $\epsilon$   small enough we have
\begin{equation}
\label{e.es.117}
\sum_j\norma{\eta_j\E_j}_{L^1_t}\leq CC_2\epsilon^2
\end{equation}
\end{lemma}
The important fact is that the right hand side is only linear in
$C_2$. The   proof of this lemma is postponed to Appendix
\ref{prle117}.

\subsection{End of the proof of Theorem \ref{proposition:mainbounds}}
\label{subsec:end}   Using the notations of section \ref{model}, for
solutions of the system \eqref{e.es.112} we have
\begin{equation}
\label{e.es.120} \frac{dH_{0L}}{dt}= -\sum_{\lambda\in\Lambda}
\langle F_{\lambda}; B_\lambda \bar F_\lambda
\rangle+\sum_j\omega_j(\eta_j \bar \E_j+\bar \eta_j \E_j)
\end{equation}
Integrating and reorganizing we get
\begin{equation}
\label{e.es.1201} H_{0L}(t)+\sum_\lambda \int_0^t\langle
F_{\lambda}; B_\lambda \bar F_\lambda \rangle(s) ds=
H_{0L}(0)+\int_0^t \sum_j\omega_j(\eta_j \bar \E_j+\bar \eta_j
\E_j)(s) ds.\nonumber
\end{equation}
Using the positivity of $H_{0L}$, we immediately get
\begin{equation}
\label{e.es.1212}
\sum_{\lambda}\int_0^t\langle F_\lambda;B_\lambda
F_\lambda\rangle(s)ds\leq (C+CC_2)\epsilon^2\ ,
\end{equation}
from which, using assumption (H7), we get
\begin{equation}
\label{e.es.121} \sum_{\mu\in M}\int_0^T |\eta^\mu|^2dt\leq
(C+CC_2)\epsilon^2, \nonumber
\end{equation}
which   implies
\begin{equation}
\label{e.es.122} \sum_{\mu\in M}\int_0^T |\xi^\mu| ^2dt\leq
(C+CC_2)\epsilon^2. \nonumber
\end{equation}
We have thus proved the following final step of the proof:
\begin{theorem}
\label{maintech}
The inequalities \eqref{4.4a} and \eqref{4.4} imply
\begin{eqnarray}
&   \|  f \| _{L^r_t([0,T],W^{ {1}/{p}- {1}/{r} ,p}_x)}\le
  K_1(C_2) \epsilon \text{ for all admissible pairs $(r,p)$} \label{4.4b}
\\& \| \xi ^\mu \| _{L^2_t([0,T])}\le
 C \sqrt{C_2} \epsilon \text{ for all multi indexes $\mu$
  with  $\omega \cdot \mu >m $} \label{4.4c}
\end{eqnarray}
\end{theorem}
Thus, provided that $C_2/2>C\sqrt{C_2} $ and $C_1/2>  K_1(C_2)$, we
see that \eqref{4.4a}--\eqref{4.4} imply the same estimates but with
  $C_1$, $C_2$ replaced by $C_1/2$, $C_2/2$. Then \eqref{4.4a} and \eqref{4.4}
hold   with $[0,T]$ replaced by $[0,\infty )$. This yields Theorem
\ref{proposition:mainbounds}.

\appendix

\section{Proofs of Lemmas \ref{lemma:weighted}
and \ref{lemma:christkiselev}}

\subsection{Proof of Lemma  \ref{lemma:weighted}}
\label{planewaves}

By a simple argument as in p.24 \cite{sofferweinsten1} which uses
Theorem \ref{yajima}, it is enough to prove, that, for any fixed $\chi
\in C^\infty _0((m,\infty ), \mathbb{R})$ with $\chi \equiv 1$ in
$[a,b]$, we have for $s>9/2$

\begin{equation}\label{inequality:bounded}
\|  \chi (B)e^{ -\im B t}R_{B}(\mu+\im 0) g\|_{L^{2,-s}_x} \le
C\langle t \rangle ^{-\frac 32} \| g \|_{L^{2,s} _x},
\end{equation}
for some fixed $C$ which depends on $\chi$.  Indeed, for
  $\overline{\chi} =1-\chi$, for any $\mu \in [a,b]$, for $s>3/2$ and
  for a fixed small $\eta >0$, there is $C$ such that, for
  $B_0=\sqrt{-\Delta +m^2}$

  \begin{equation}\label{inequality:bounded1} \begin{aligned}&
\| \overline{ \chi} (B)e^{ -\im B t}R_{B}(\mu ) g\|_{H^{-4,-s}_x} \le  \| \overline{ \chi} (B)e^{ -\im B t}R_{B}(\mu ) g\|_{W^{-4+\eta, \infty}_x  }\\&
\le  C_3\| \overline{ \chi} (B_0)e^{ -\im B_0 t}R_{B_0}(\mu )\mathcal{Z} _+g\|_{W^{-4+\eta, \infty}_x  }\le
 C_2  \langle t \rangle ^{-\frac 32}
 \| \mathcal{Z} _+   g\|_{L^{1}_x  }
\\& \le
 C_1  \langle t \rangle ^{-\frac 32}
 \|    g\|_{L^{1}_x  }\le
C\langle t \rangle ^{-\frac 32} \| g \|_{L^{2,s} _x}, \end{aligned}
\end{equation}
  for all $g\in
  L^{2,s} _x$. So we focus on     \eqref{inequality:bounded}.   We have

 \begin{equation}\label{term1} \begin{aligned}
&\langle x \rangle ^{-\gamma } \chi (B) e^{-\im Bt}R ^+( \mu )
 \langle y \rangle ^{-\gamma }= \\& \lim _{\epsilon \searrow 0}
   e^{-\im \mu t}\langle x
\rangle ^{-\gamma } \int _t^{+\infty }e^{-\im (B-\mu -\im \epsilon )s}
 \chi (B)
  ds\langle y \rangle ^{-\gamma }.
\end{aligned}
\end{equation}
Using the distorted plane waves $u(x,\xi )$ associated to the
continuous spectrum of $-\Delta +V$, we can write the following integral kernel:
\begin{equation}\label{term2} \begin{aligned}&\langle x \rangle ^{-\gamma } \left ( \chi (B)
e^{-\im (B-\mu -\im \epsilon )s}\right ) (x,y) \langle y \rangle ^{-\gamma }=
  \\&  \langle x \rangle ^{-\gamma } \int _{\Bbb
R^3} u(x,\xi )e^{(-\im  \sqrt{\xi ^2 +m^2} +\im  \mu -
\epsilon)s}\chi (\sqrt{\xi ^2 +m^2} ) \bar u(y,\xi ) d\xi \langle y
\rangle ^{-\gamma }.
\end{aligned}
\end{equation}
We have $ u(x,\xi ) = e^{\im x\cdot \xi }   + e^{\im x\cdot \xi
}w(x, \xi ),$ with $w(x, \xi )$ the unique solution in $L^{2,-s}$,
$s> 1/2$, of the integral equation
\begin{equation}\label{eq:lippmanschwinger1}  w(x,\xi )= -F(x,\xi )-
\int _{\Bbb R^3}w(y, \xi )V  (y)
  \frac{e^{\im |\xi | |y-x|}}{4\pi |y-x|}
e^{\im (y-x)\cdot \xi } dy,
\end{equation}
with
\begin{equation}\label{eq:lippmanschwinger2}F(x,\xi )=
\int _{\Bbb R^3}V(y)
 \frac{e^{\im |\xi | |y-x|}}{4\pi |y-x|} e^{i(y-x)\cdot \xi } dy.
\end{equation}   It is
elementary to show that $|V(x)|\le C \langle x \rangle ^{-5-\sigma}
$ for $\sigma >0$ implies that, for $\xi $ in the support of $\chi
(\sqrt{\xi ^2 +m^2} )$ and for $|\alpha |\le 3$,  then $|\partial
^\alpha _ \xi F(x,\xi )| \le \tilde c_{\alpha   } \langle x \rangle
^{|\alpha |-1} $ for fixed constants $\tilde c_{\alpha   }$. By
elementary arguments, as in \cite{cu2}, from stationary scattering
theory it is possible for $|\alpha |\le 3$ to conclude
correspondingly $|\partial ^\alpha _\xi w(x,\xi )| \le c_{\alpha }
\langle x \rangle ^{|\alpha |-1} $ for fixed constants $c_{\alpha }$.
Then, using  $e^{ -\im s
\sqrt{\xi ^2 +m^2}}= \frac{\im \sqrt{\xi ^2 +m^2} }{|\xi | s} \frac
d {d|\xi |} e^{ -\im s \sqrt{\xi ^2 +m^2} }$ we have

 \begin{equation} \begin{aligned}&
 \text{rhs\eqref{term2}}= (-1)^r \langle x \rangle ^{-\gamma }  \langle y
\rangle ^{-\gamma } \times \\&  \int _{ \R ^3}e^{(-\im  \sqrt{\xi ^2 +m^2} +\im  \mu -
\epsilon)s}       \left ( \frac \partial {\partial |\xi |}   \frac{\im \sqrt{\xi ^2 +m^2} }{|\xi | s} \right ) ^r \left [ u(x,\xi )\chi (\sqrt{\xi ^2 +m^2} ) \bar u(y,\xi ) \right ] d\xi .\end{aligned}
\nonumber
\end{equation}
This yields
 \begin{equation} \begin{aligned} &|\text{rhs\eqref{term2}} |\le c\langle x \rangle ^{-\gamma +r}
\langle y \rangle ^{-\gamma +r} s^{-r}e^{-\epsilon t}\, \text{ and
so} \\& |\text{rhs\eqref{term1}}|\le c\langle x \rangle ^{-\gamma +r}
\langle y \rangle ^{-\gamma +r} t^{-r+1}.\end{aligned}\nonumber
\end{equation} For $\gamma > r+3/2$ and
$r=3$, we obtain the conclusion.

\begin{remark}
\label{remark:planewaves} Notice that when $|V(y)|\le C e^{-a|y|}$
for $a>0$, equations
\eqref{eq:lippmanschwinger1}--\eqref{eq:lippmanschwinger2} make
sense with $\im |\xi |$ replaced by $\sqrt{  - \xi _1^2- \xi _2^2-
\xi _3^2}$ with $\xi $ in an open neighborhood $U$ of $\R
^3\backslash \{ 0\}$  in $\C ^3\backslash \{ 0\}$. Then we get
solutions $w(x,\xi )$ bounded and analytic in $\xi$. Correspondingly
we obtain   $u(x,\xi )$ for $\xi \in U$ analytic in $U$ and with
$|u(x,\xi )|\le C e^{|x|\sum _{j=1}^{3}|\Im \xi _j|}$. Consequently,
if $|v(x)|\le c_0e^{-b|x|}$ for $b>0$ and for the distorted plane
wave transformation
\begin{equation}\label{distortedfourier}
\widehat{v}(\xi )= (2\pi )^{-\frac{3}{2}}\int _{\R ^3}\overline{u}
(y,\xi ) v(y) dy,
\end{equation}
then $\widehat{v}(\xi )$ extends into an holomorphic function in
some open neighborhood   of $\R ^3\backslash \{ 0\}$  in $\C
^3\backslash \{ 0\}$.

\end{remark}

\subsection{Proof of Lemma  \ref{lemma:christkiselev}}

The proof originates from \cite{mizumachi} (in fact see also
\cite{RoSc}) but here we state the steps of a simplification in
\cite{cuccagnatarulli}. We first state Lemmas
\ref{lemma:Katosmoothness1}--\ref{lemma:Katosmoothness2}. They imply
Lemma  \ref{lemma:christkiselev} by an argument in \cite{mizumachi}.
First of all we need some estimates on the resolvent, for the proof
see Lemma 2.8 \cite{danconafanelli}:

\begin{lemma}
\label{lemma:Katosmoothness1}  For any $s>1$  there is a $C>0$ such
that for any $z$ with $\Im z>0$ we have
\begin{equation}\label{eq:Katosmoothness1}
\|  R_{B}(z)P_c \|_{B(L^{2,s} _x, L^{2,-s}_x)} \le C .
\end{equation}
\end{lemma}
  Estimates \eqref{eq:Katosmoothness1} yield a  Kato smoothness
  \cite{kato} result,
  see the proof of Lemma 3.3 \cite{cuccagnatarulli}:

\begin{lemma}
\label{lemma:Katosmoothness2} For any $s>1$ there is a $C$ such that for
all Schwartz functions $u_0(x)$ and $g(t,x)$ we have
\begin{eqnarray}\label{eq:Katosmoothness2} &
\|  e^{-\im Bt}P_c u_0\|_{ L^2_t L^{2,-s}_x} \le C \| P_c u_0\|_{
L^{2}_x}\\ \label{eq:Katosmoothness3}   & \left\|\int_\Bbb R e^{\im
tB }P_c g(t,\cdot)dt\right\|_{L^2_x } \le C  \| P_c g\|_{ L_t^2
L^{2, s}_x} .
\end{eqnarray}
\end{lemma}
Now we are ready to prove Lemma \ref{lemma:christkiselev}.
  For $g(t,x
)\in C^\infty _0(  \mathbb{R}\times  \mathbb{R}^3)$ set
\begin{eqnarray}T g(t)=\int _0^{+\infty}
e^{-\im (t-s)B }P_c g(s) ds  .\nonumber \end{eqnarray}
\eqref{eq:Katosmoothness3} implies $f:=\int _0^{+\infty}  e^{\im sB
}P_c g(s)ds\in L^2_x$.  Lemma \ref{lemma-strichartz}   implies that
for all $(p,q) $ admissible we have
\begin{eqnarray}  \|T g(t)\|_{L ^{  p}_t
W_x ^{ \frac{1}{q}-\frac{1}{p} , q}  }\lesssim  \| f\|_{H
^{\frac{1}{2}}_x} \lesssim \| \sqrt{B} P_c g\|_{L^2_tL_x ^{2, s}  }
 \nonumber
\end{eqnarray}
where the last inequality follows from \eqref{eq:Katosmoothness3}
and Theorem \ref{yajima}:
\begin{equation}  \begin{aligned} &  \| f\|_{H
^{\frac{1}{2}}_x} \lesssim \| \mathcal{Z} _+ f\|_{H
^{\frac{1}{2}}_x} \lesssim \| (-\Delta
+m^2)^{\frac{1}{4}}\mathcal{Z} _+ f\|_{L ^{2}_x}
\\& \lesssim \|  \sqrt{B} f\|_{L ^{2}_x} \lesssim
\| \sqrt{B} P_c g\|_{L^2_tL_x ^{2, s}  }, \text{ by $\sqrt{B}f =\int _0^{+\infty}  e^{\im sB }\sqrt{B}P_c g(s)ds\in
L^2_x$.}
\end{aligned} \nonumber
\end{equation}
 Notice that \eqref{eq:Katosmoothness3} implies also $\|
f\|_{H ^{\frac{1}{2}}_x} \lesssim \| \sqrt{B} P_cg\|_{L^a_tL_x ^{2,
s} } $ for any $a\in [1,2]$. The following
  well known  results by Christ \& Kieselev,
  see  Lemma 3.1 \cite{soggesmith},
  yields Lemma \ref{lemma:christkiselev}.
\begin{lemma}
Consider two Banach spaces and $X$ and $Y$ and $K(s,t)$ continuous
function valued in the space $B(X, Y)$. Let
$$T_K f(t)=\int_{-\infty}^\infty K(t,s)f(s) ds\text{ and }
 \tilde T_K f(t) = \int_{-\infty}^t K(t,s) f(s) ds.$$
 Then we have:
 Let $1\leq a<  b \leq \infty$ and $I$ an interval.
Assume that there exists $C>0$ such that
$$
\|T_K f\|_{  L^b(I,Y) }\leq C \|f\|_{L^a(I,X)}.$$ Then
$$
\|\tilde T_K f\|_{  L^b(I,Y) }\leq C' \|f\|_{L^a(I,X)} $$ where
$C'=C'(C, a,b )>0$.
\end{lemma}

\section{Proof of Lemma \ref{l.es.11}. }\label{prle117}

First of all   \eqref{4.4} immediately implies the   estimate
\begin{equation}
\label{e.dis.1}
 \| \eta ^\mu \| _{L^2_t([0,T])}\le
2 C_2 \epsilon \text{ for all multi indexes $\mu$
  with  $\omega \cdot \mu >m $.}
\end{equation}
Let us start with the contribution of the the last line of
\eqref{e.es.bis}.

\begin{lemma}
\label{lemma:dimenticato} We have
\begin{equation}
\label{eq:dimenticata} \norma{ \left [- \frac{\partial
Z_0}{\partial\bar\xi_j}(\xi)+ \N_j
    (\xi)+ \frac{\partial Z_0}{\partial\bar\xi_j}(\eta)- \N_j (\eta)
\right ] \eta _j }_{L^1_t}\leq C  \epsilon^3.
\end{equation}\end{lemma}
\proof For definiteness we focus on    $\| (\overline{\partial}
_jZ_0(\xi   ) -\overline{\partial} _jZ_0(\eta   ))\bar \eta _j\|
_{L^1_t}$.  It is enough to consider quantities  $  \xi ^{\alpha}
\frac{\bar \xi ^{ {\beta} }}{\bar \xi _j} \bar \eta _j   - \eta
^{\alpha} \frac{\bar \eta ^{ {\beta} }}{\bar \eta _j} \bar \eta _j   $ with $\omega \cdot \alpha =
\omega \cdot \beta  $ and $\beta _j>0$. By Taylor expansion these
are
\begin{equation}  \sum  _k\partial _k
\left (\frac{\xi  ^{\alpha} \bar \xi  ^{ {\beta} }}{\bar \xi
_j}\right ) (\eta _k-\xi  _k)\bar \eta _j + \sum
_k\overline{\partial} _k \left (\frac{\xi ^{\alpha} \bar \xi ^{
{\beta} }}{\bar \xi _j}\right ) (\bar \eta  _k-\bar \xi _k)\bar \eta
_j +\bar \eta _j O(|\xi -\eta |^2).\nonumber
\end{equation}
The reminder term is the easiest, the other two terms similar.
Substituting  \eqref{equation:FGR2}, a typical term in the first
summation is $ \frac{\xi ^{\alpha +A} \bar \xi ^{B+\beta}    }{|\xi _k|^2} $, with  all four $\alpha$, $\beta$, $A$ and
$B$ in $M$ and with $\alpha _k\neq 0\neq B_k$.  (H5) and $\omega
\cdot \alpha =\omega \cdot \beta$ imply that
   there is at least
one index $\beta _\ell \neq 0$ such that $\omega _\ell = \omega _k$.
Then
\begin{equation}\label{equation:FGR12} \left \| \frac{\xi ^{\alpha} \bar \xi ^{\beta}  \xi ^A \bar
\xi ^{B} }{|\xi _k|^2}\right \| _{L^1_t}\le \left \| \xi ^{ A}
\right  \| _{L^2_t} \left \| \frac{ \xi ^{  {B} }  \xi _\ell}{\xi
_k} \right \| _{L^2_t} \left \| \frac{\xi ^{ {\alpha}  } \bar \xi ^{
{\beta}  }}{\xi _\ell \xi _k}\right \| _{L^\infty_t}\lesssim
C_2^2\epsilon ^{ |\alpha|+|\beta |}\le C_2^2\epsilon ^{4}
\end{equation}
by the fact that   monomials $\xi ^\alpha \bar \xi ^\beta $ in $Z_0$
are such that $|\alpha |=|\beta |\ge 2$. Other terms can be bounded
similarly.\qed

\begin{lemma}
\label{l.dis.s} For $\epsilon$   small enough  we have
\begin{equation}
 \norma{  \eta _j    \langle  \partial _{  \overline{\xi}_j}
G
  ,g \rangle
}_{L^1_t}   + \norma{  \eta _j    \langle  \partial _{
\overline{\xi}_j}   \bar G
  ,\overline{g} \rangle
}_{L^1_t}   \leq C C_2\epsilon^2. \nonumber
\end{equation}
\end{lemma}

\proof We  first bound $ \norma{  \eta _j    \langle  \partial
_{\bar \xi_j} G
  ,g \rangle
}_{L^1_t}$.  We have by Lemma \ref{lemma:bound g}

$$\norma{  \eta  _j  \left\langle  \partial
_{\bar \xi_j} G
  ,g\right\rangle
}_{L^1_t}\le \norma{\eta _j  \partial _{\bar \xi_j}
G}_{L^2_tH^{4,s}} \norma{g}_{L^2_tH^{-4,-s}}\le
   C_0\epsilon \norma{\eta _j  \partial _{\bar
\xi_j}   G}_{L^2_tH^{4,s}}.$$ We have \begin{equation} \label{scomp
eta1} \norma{\eta _j   \partial _{\bar \xi_j} G}_{L^2_tH^{4,s}}\le
\norma{\xi _j   \partial _{\bar \xi_j} G}_{L^2_tH^{4,s}}+
\norma{\Delta _j   \partial _{\bar \xi_j} G}_{L^2_tH^{4,s}}
.\end{equation}  By \eqref{m.2}--\eqref{m.3} and \eqref{m.17bis} we
have
\begin{equation} \label{scomp
eta2}  \begin{aligned} &  \norma{\Delta _j   \partial _{\bar \xi_j}
G}_{L^2_tH^{4,s}} \le  \norma{\Delta _j   }_{L^2_t } \norma{
\partial _{\bar \xi_j} G}_{L^\infty_tH^{4,s}}\\&
\le C \sum _{\mu \in M} \norma{\xi ^\mu }_{L^2_t } \norma{\xi
}_{L^\infty _t }^2\le C  C_2  \epsilon
^3.\end{aligned}\end{equation} Finally,  by \eqref{m.2}--\eqref{m.3}
 we have

\begin{equation}  \norma{\xi _j   \partial _{\bar \xi_j}
G}_{L^2_tH^{4,s}} \leq      C \sum_{\omega \cdot (\nu  -\mu
)>m}\norma{\xi^\mu\bar\xi^\nu   }_{L^2_t }  \leq CC_2\epsilon
.\nonumber
\end{equation}
Now we bound $ \norma{  \eta _j    \langle  \partial _{\bar \xi_j}
\bar {G}
  ,\overline{g} \rangle
}_{L^1_t}$. We reduce to an analogue of \eqref{scomp
eta1}--\eqref{scomp eta2}
\begin{equation}   \begin{aligned} &\norma{\eta _j   \partial _{\bar \xi_j}\bar
G}_{L^2_tH^{4,s}}\le \norma{\xi _j   \partial _{\bar \xi_j} \bar
G}_{L^2_tH^{4,s}}+ \norma{\Delta _j   \partial _{\bar \xi_j}\bar
G}_{L^2_tH^{4,s}} \\& \le \norma{\xi _j   \partial _{  \xi_j}
G}_{L^2_tH^{4,s}}+ CC_2\epsilon ^2.\end{aligned}\nonumber
\end{equation} Finally
\begin{equation} \begin{aligned} &\norma{\xi _j   \partial _{  \xi_j}
G}_{L^2_tH^{4,s}} \lesssim \sum _{\omega \cdot \nu >m } \norma{\mu
_j   \xi ^\mu \bar \xi ^\nu}_{L^2_t }\le CC_2 \epsilon ^2.
\end{aligned} \nonumber
\end{equation} \qed

\begin{lemma}
\label{Delta} For $\epsilon$   small enough  we have
\begin{equation}
  \norma{  \eta _j     \partial _{\xi _k} \Delta _j }_{L^2_t} +
  \norma{  \eta _j     \partial _{\overline{\xi} _k}
  \Delta _j }_{L^2_t}     \leq C
C_2\epsilon^2. \nonumber
\end{equation}
\end{lemma}
\proof We first bound $\norma{  \eta _j     \partial _{\xi _k}
\Delta _j }_{L^2_t} $. As in \eqref{scomp eta1}--\eqref{scomp eta2}
we write
\begin{equation}   \norma{\eta _j   \partial _{\xi _k} \Delta _j }_{L^2_t}\le
\norma{\xi _j    \partial _{\xi _k} \Delta _j }_{L^2_t}+
\norma{\Delta _j   \partial _{\xi _k} \Delta _j }_{L^2_t} \le
\norma{\xi _j    \partial _{\xi _k} \Delta _j }_{L^2_t} +C C_2
\epsilon ^2 .\nonumber\end{equation} We have
$$\xi_j\frac{\partial \Delta_j}{\partial
  \xi_k} \sim
\frac{ \xi^\mu\bar\xi^{\nu}}{\xi_k}   \text{ with $\mu    ,\nu$ in $
M$,
  $\mu _k\neq 0$ }.
 $$
Then, by $\mu _k\neq 0$ and $|\mu|\ge 2$,  we have
\begin{equation}
\label{e.dis.7}\left \| \frac{ \xi^\mu\bar\xi^{\nu}}{\xi_k} \right
\| _{L^2_t}\le \|  \xi^{\nu}   \| _{L^2_t}\left \| \frac{ \xi^{\mu }
}{\xi_k}
  \right \| _{L^\infty _t}     \le CC_2 \epsilon ^2.\end{equation}
Now we bound $\norma{  \eta _j     \partial _{\overline{\xi} _k}
  \Delta _j }_{L^2_t}\le \norma{  \xi _j     \partial _{\overline{\xi} _k}
  \Delta _j }_{L^2_t} +C C_2\epsilon ^2$. We have
$$\xi_j\frac{\partial \Delta_j}{\partial
  \bar \xi_k} \sim
\frac{   \xi^\mu\bar\xi^{\nu}}{\bar \xi_k}   \text{ with $\mu
,\nu$ in $ M$,
  $\nu _k\neq 0$ }.
 $$
We then exploit
\begin{equation}
 \left \| \frac{   \xi^\mu\bar\xi^{\nu}}{\bar \xi_k}  \right
\| _{L^2_t}\le \|  \xi^{\mu}   \| _{L^2_t}\left \|  \frac{{\xi ^\nu}}{\xi_k}
  \right \| _{L^\infty _t}     \le CC_2 \epsilon ^2.\nonumber \end{equation}

\qed
\begin{lemma} \label{G1k} We have
$  \norma{  \eta _j   \mathcal{G}_{1, {j}} }_{L^1_t}\leq C
(C_2)\epsilon^3 $  and $  \norma{      \mathcal{G}_{1, {j}}
}_{L^2_t}\leq C (C_2)\epsilon^2 $.
\end{lemma}
\proof As in \eqref{scomp eta1}  we write
\begin{equation}   \norma{\eta _j   \mathcal{G}_{1, {j}} }_{L^1_t}\le
\norma{\xi _j    \mathcal{G}_{1, {j}} }_{L^1_t}+ \norma{\Delta _j
\mathcal{G}_{1, {j}} }_{L^1_t} .\nonumber\end{equation}
 $ |\xi_j\mathcal{G}_{1,j} |$  is bounded by  the absolute values of terms of
the form either
\begin{equation}
\label{case1} \xi^{\mu+\mu'}\bar \xi^{ \nu +\nu '}\, , \, \mu '\in M
\, , \, \nu \in M\, , \,   (\mu , \nu ')\neq (0,0),
\end{equation}
which originate from terms in \eqref{m.14bis} with $(\mu , \nu ')\neq (0,0)$,
or by terms originating  from terms in \eqref{m.14},
\begin{equation}
\label{case2} \xi _j \xi^{ \mu'}\bar \xi^{ \nu  }\, , \, \mu '\in M
\, , \, \nu \in M.
\end{equation} In case \eqref{case1}

$$
\| \xi^{\mu+\mu'}\bar \xi^{ \nu +\nu '}\| _{L^1_t}\le \| \xi^{\nu
}\| _{L^2_t} \|   \xi^{\mu' }\| _{L^2_t}\| \xi \| _{L^\infty
_t}^{|\mu |+|\nu '|}\le C  C_2^2\epsilon ^3.
$$
Similarly, in case \eqref{case2}

$$
\| \xi _j \xi^{ \mu'}\bar \xi^{ \nu  }\| _{L^1_t}\le \| \xi^{\nu }\|
_{L^2_t} \|   \xi^{\mu' }\| _{L^2_t}\| \xi _j \| _{L^\infty _t} \le
C  C_2^2\epsilon ^3.
$$
Dividing \eqref{case1}--\eqref{case2} by $\xi _j$ we see that
 \begin{equation}   \norma{  \mathcal{G}_{1, {j}}
}_{L^2_t}\le C  C_2 \epsilon ^2.\nonumber
\end{equation}
Finally, $\norma{\Delta _j \mathcal{G}_{1, {j}} }_{L^1_t}\le
\norma{\Delta _j   }_{L^2_t}\norma{  \mathcal{G}_{1, {j}}
}_{L^2_t}\le C  C_2^2\epsilon ^3.$

 \qed

\begin{lemma}
\label{G 0} We have  $  \norma{ \mathcal{G}_{0,k}(\xi )
}_{L^2_t}\leq C  C_2 \epsilon^2 .$
\end{lemma}
\proof Indeed by \eqref{m.13ter}, \eqref{4.4} and remark
\ref{r.es.2} we have
\begin{equation} \norma{
\mathcal{G}_{0,j}(\xi ) }_{L^2_t}\le \sum _{  \mu , \nu \in M  } \nu
_j \norma{  \frac{\xi ^\mu \xi ^\nu}{\xi _j} }_{L^2_t} \le C
C_2\epsilon^2.\nonumber
\end{equation}
\qed

\begin{lemma}
We have: \label{l.es.333}
\begin{equation}
  \norma{\eta _j (\partial _{ \xi_l} \Delta_j)
(\partial _{  \bar \xi_l}
  Z_0) }_{L^1_t}\leq CC_2^2\epsilon^3\ .\nonumber
\end{equation}\end{lemma}
\proof We have
\begin{equation} \label{e.es.334}
\norma{\eta _j (\partial _{ \xi_l} \Delta_j)
(\partial _{  \bar \xi_l}
  Z_0) }_{L^1_t} \le \norma{\xi _j (\partial _{ \xi_l} \Delta_j)
(\partial _{  \bar \xi_l}
  Z_0) }_{L^1_t} +\norma{\Delta _j (\partial _{ \xi_l} \Delta_j)
(\partial _{  \bar \xi_l}
  Z_0) }_{L^1_t}.\end{equation}
We first bound the first term in rhs of \eqref{e.es.334}. It has a sum
of terms of the form
\begin{equation}
\label{e.es.335}  \frac{\xi^\mu\bar
\xi^{\nu}}{\xi_l}\frac{\xi^\alpha\bar
  \xi^\beta}{\bar \xi_l}
\end{equation}
with indexes   such that
\begin{equation}
\label{e.es.336}
 \text{$\mu$ and $  \nu\in M$,
$\omega\cdot(\alpha-\beta)=0$, $\mu _l\neq 0 \neq \beta _l$.}
\end{equation}
By (H5) there is $\alpha _k\neq 0$ such that $\omega _k =\omega _l$.
Then
\begin{equation}
\label{e.es.339}\left \|  \frac{\xi^\mu\bar
\xi^{\nu}}{\xi_l}\frac{\xi^\alpha\bar
  \xi^\beta}{\bar \xi_l}\right \| _{L^1_t} \le
  \left \|  \xi ^{\nu}\right \| _{L^2_t}
  \left \|  \frac{\xi ^{\mu  }\xi _k}{\xi _l}\right \| _{L^2_t}
  \left \|   \frac{\xi^\alpha\bar
  \xi^\beta}{\xi _k \bar \xi_l}\right \| _{L^\infty _t}\lesssim
  C_2^2\epsilon ^{|\alpha|+|\beta|}\le C_2^2\epsilon ^{4}
\end{equation}
by the fact that   monomials $\xi ^\alpha \bar \xi ^\beta $ in $Z_0$
are such that $|\alpha |=|\beta |\ge 2$.

Finally, by  \eqref{m.17bis}

\begin{equation}  \begin{aligned} &
\norma{\Delta _j (\partial _{ \xi_l} \Delta_j)
(\partial _{  \bar \xi_l}
  Z_0) }_{L^1_t}\le \norma{\Delta _j}_{L^2_t}
  \norma{(\partial _{ \xi_l} \Delta_j)
(\partial _{  \bar \xi_l}
  Z_0) }_{L^2_t}\\& \le C C_2
  \epsilon ^2\norma{ (\partial _{ \xi_l} \Delta_j)
(\partial _{  \bar \xi_l}
  Z_0) }_{L^2_t}.
\end{aligned} \nonumber
\end{equation}
The last factor can be bounded using
\begin{equation} \label{e.es.400} \begin{aligned} & \norma{\frac{\xi^\mu\bar
\xi^{\nu}}{\bar \xi _j\xi_l}\frac{\xi^\alpha\bar
  \xi^\beta}{\bar \xi_l} }_{L^2_t} \le
  \norma{\frac{\xi^\mu \xi _k}{ \xi_l}  }_{L^2_t}
  \norma{\frac{\xi^\alpha \bar
\xi^{\nu +\beta }}{\bar \xi _j \bar \xi_l \xi _k }  }_{L^\infty _t}
   \le CC_2\epsilon ^{2},
\end{aligned}
\end{equation}
where in the last formula the exponents satisfy \eqref{e.es.336} and
$\mu _j\neq 0$ and where we picked $k$ such that  $\alpha _k\neq 0$ and  $\omega _k =\omega _l$. \qed

\begin{lemma}
We have: \label{l.es.334}
\begin{equation}
  \norma{\eta _j (\partial _{ \bar \xi_l} \Delta_j)(
\partial _{    \xi_l}
  Z_0 )}_{L^1_t}\leq CC_2^2\epsilon^3\ .\nonumber
\end{equation}\end{lemma}
\proof We have
\begin{equation}
\norma{\eta _j (\partial _{ \bar \xi_l} \Delta_j)(
\partial _{    \xi_l}
  Z_0 ) }_{L^1_t} \le \norma{\xi _j (\partial _{ \bar \xi_l} \Delta_j)(
\partial _{    \xi_l}
  Z_0 ) }_{L^1_t} +\norma{\Delta _j (\partial _{ \bar \xi_l} \Delta_j)(
\partial _{    \xi_l}
  Z_0 )  }_{L^1_t}.\nonumber \end{equation}
We first bound the first term in rhs. It has a sum of terms of the
form
\begin{equation}
\label{e.es.337}  \frac{\xi^\mu\bar \xi^{\nu}}{\bar
\xi_l}\frac{\xi^\alpha\bar
  \xi^\beta}{  \xi_l}
\end{equation}
with indexes   such that
\begin{equation}
\label{e.es.338}
 \text{$\mu$ and $  \nu\in M$,
$\omega\cdot(\alpha-\beta)=0$, $\nu _l\neq 0 \neq \alpha _l$.}
\end{equation}
Since  complex conjugates  of terms
\eqref{e.es.337}--\eqref{e.es.338} give terms
\eqref{e.es.335}--\eqref{e.es.336}, we get the desired estimates  by
\eqref{e.es.339}. By this argument and by   \eqref{e.es.400},  $\norma{\partial _{ \bar \xi_l}
\Delta_j
\partial _{    \xi_l}
  Z_0 }_{L^2_t}\le CC_2 \epsilon ^2 $. \qed

\medskip
Finally we complete the proof of Lemma \ref{l.es.11}. The
contribution from the last  line of \eqref{e.es.bis} is bounded in
Lemma \ref{lemma:dimenticato}. We have    $\| \partial _{\bar \xi}
\resto \| _{L^1_t}\le C(C_1) \epsilon ^2$, see  \eqref{e.es.6}. Then
$\| \eta \partial _{\bar \xi} \resto \| _{L^1_t}\le C(C_1) \epsilon
^3< c_0  \epsilon ^2$ for any preassigned $c_0$. Hence all the
related terms in $\eta _j \mathcal{E} _j$ satisfy a better estimate
than   \eqref{e.es.117}. Similarly,

\begin{equation} \| \eta \partial _{  \xi} \resto \| _{L^1_t} =
\| \overline{\eta} \partial _{  \overline{\xi}} \resto \| _{L^1_t}
 =
\|  {\eta} \partial _{  \overline{\xi}} \resto \| _{L^1_t}\le C(C_1)
\epsilon ^3< c_0  \epsilon ^2 .\nonumber
\end{equation}
So we are left with  the contributions of the   terms  in the first
three lines in \eqref{e.es.bis} not coming from   $   \resto$.

The contribution from $( \partial _{ \xi_l} \Delta_j)(
\partial _{    \bar \xi_l}
  Z_0) $ in the first line of   \eqref{e.es.bis} is bounded in
Lemma \ref{l.es.333}.  The other terms from the first line of
\eqref{e.es.bis} are bounded  by
\begin{equation}  \begin{aligned} & \| \eta _j(\partial _{
\xi_k} \Delta_j )\langle \partial _{ \bar \xi_k} G, f \rangle \|
_{L^2_t}+\| \eta _j(\partial _{
\xi_k} \Delta_j ) \langle \partial _{
\bar \xi_k} \bar G, \bar f \rangle \| _{L^1_t}\\& \le \| \eta
_j\partial _{ \xi_k} \Delta_j\| _{L^2_t} \| f \|
_{L^2_tW_x^{-1/3,6}} \le C C_2 C_1\epsilon ^3,
\end{aligned}\nonumber
\end{equation}
where we have used \eqref{4.4} and Lemma \ref{Delta}. In Lemmas \ref{l.dis.s} and \ref{G1k} we have bounded the
contributions from the second line of \ref{e.es.bis} coming from the
$\delta _{jk}$, and all terms in   $\mathcal{G}_{1,k}$. The
remaining terms, thanks to Lemma \ref{lemma:bound g}   are bounded
by
\begin{equation}  \begin{aligned} & \| \eta _j(\partial _{
\xi_k} \Delta_j) \langle \partial _{ \bar \xi_k} G, g \rangle \|
_{L^2_t}+\| \eta _j(\partial _{
\xi_k} \Delta_j) \langle \partial _{
\bar \xi_k} \bar G, \bar g \rangle \| _{L^1_t}\\& \le \| \eta
_j\partial _{ \xi_k} \Delta_j\| _{L^2_t} \| g \| _{L^2_tH_x^{-4,-s}}
\le C C_2  \epsilon ^3.
\end{aligned}\nonumber
\end{equation}
Focusing on the third line of \eqref{e.es.bis}, the terms from
$\partial _{ \bar \xi_k} \Delta_j  \partial _{ \xi_k} Z_0   $ are
bounded by Lemma \ref{l.es.334}. The other terms, by Lemmas
\ref{Delta}, \ref{G1k}, \ref{G 0}.

\section{Regularization estimates and proof of Lemma \ref{sol.homo1}. }\label{app.resolv}

First of all Lemma \ref{lem:m.11} is a consequence of the following lemma.
\begin{lemma}\label{lemma-regularization}
Let $|V(x)|\le C \langle x \rangle ^{-5}$. Then, for $\Phi \in H^{2,
s}$ for $s>1/2$ and    $\lambda >m$, we have that   $ R^{\pm}_B(\lambda
)\Phi $ are well defined and belong  to $L^{2, -s}$.
\end{lemma}

\proof  We   set $\Psi=(B+\lambda )\Phi $. Then $Y=R^{+}_B(\lambda)
\Phi = R^{+}_{-\Delta +V}(k^2) \Psi $ with $k^2=\lambda ^2-m^2$ (the
proof for $R^{-}_B(\lambda) \Phi $ is similar). $|V(x)|\le C \langle
x \rangle ^{-5}$ implies that $V(x)$ is an Agmon potential, see
Example 2 XIII.8  \cite{reedsimon}.  So if $\Psi \in L^{2, s}$ then
$R^{+}_{-\Delta +V}(k^2) \Psi $ is well defined and in $L^{2, -s}$,
see Theorem XIII.33 \cite{reedsimon}. Since $\Psi \in L^{2, s}$ if  $B\Phi
\in L^{2, s}$, and since the latter is guaranteed by  Lemma \ref{lem:regularization1} below, Lemma \ref{lemma-regularization} is proved.
\qed

\begin{lemma}\label{lem:regularization1}
Let $|V(x)|\le C \langle x \rangle ^{-5}$. Then, for $\Phi \in H^{2,
s}$ for $s\ge 0 $  and for any $\kappa \in [0,1]$ we have $ B^{2\kappa}\Phi \in L^{2,  s}$.
\end{lemma}

\proof  Notice that the case $B^0=P_c$ and $B^2=(-\Delta +V)P_c$, is elementary. So we consider $\kappa \in (0,1)$.
By the  Spectral
Theorem,   for any fixed $a >0$ we write
\begin{equation}
\label{Katosquareroot} \begin{aligned} &B^{2\kappa}\Phi =c_\kappa \int _0^a
  (B^2+\tau )^{-1} B^2\Phi   \frac{d\tau
}{ \tau ^{1-\kappa }}+ c_\kappa \int _a^\infty (B^2+\tau )^{-1} B^2\Phi
\frac{d\tau }{\tau ^{1-\kappa }} , \\&
\text{with } c_\kappa = \int _0^\infty \tau ^{ \kappa -1 } (\tau +1)^{-1}d\tau .
\end{aligned}  \end{equation}
Set $
B^{2\kappa}\Phi (x)=\int _{\mathbb{R}^3} (K_a(x,y) +H_a(x,y))(B^2\Phi )(y)dy$,
with the integral kernels written in the order of the operators in
\eqref{Katosquareroot}. Set $\mathcal{H}=-\Delta +V$. We have
$B^2\Phi =(\mathcal{H}+m^2)P_c\Phi
 \in L^{2, s}$.  It is not restrictive to assume $P_c\Phi =\Phi$.
We choose $a\ge 0$ such that $V(x)+m^2+a \ge 0$ for all $x\in \R^3$
exploiting the fact that $V\in L^\infty (\R ^3)$ by (H1). Then by
the Trotter    formula, see  Theorem A.1  p.381 \cite{taylor}, we
have $e^{-t(\mathcal{H} +m^2 +\tau)} (x,y)\le e^{-t(-\Delta  +m^2
+\tau -a)} (x,y)$ for $\tau \ge a$. Then, for $\sigma =\tau -a\ge 0$
\begin{equation}
\label{Laplace1} \begin{aligned} &0<(\mathcal{H} +m^2 +\tau
)^{-1}(x,y)  = \int _0^\infty e^{-t(\mathcal{H} +m^2 +\tau)} (x,y)
dt\\& \le \int _0^\infty e^{-t(-\Delta  +m^2 +\sigma)} (x,y) dt
=(-\Delta  +m^2 +\sigma )^{-1}(x,y)\\& = \frac{e^{- \sqrt{\sigma
+m^2} |x-y|}}{4\pi |x-y|}.
\end{aligned} \nonumber \end{equation}
 Then for some fixed  constant $C>0$
\begin{equation}
\label{boundKernel}  |H_a(x,y)|\le   \int _0^\infty \frac{e^{-
\sqrt{\sigma +m^2} |x-y|}}{4\pi ^2 |x-y|}   \frac{d\sigma
}{\sigma ^{1-\kappa }} \le C  \frac{e^{-  m |x-y|/2}}{ |x-y| ^2} .
\end{equation}
By \eqref{boundKernel} we obtain that $T_s(x,y):=\langle x \rangle
^{s}\langle y \rangle ^{-s}|H_a(x,y)|$ is for any $s$  the kernel of
an operator bounded in $L^2$ by the fact that Young inequality
holds:
$$ \sup _x \| T_s(x,y) \| _{L^1_y}+\sup _y \| T_s(x,y) \| _{L^1_x}
<C_s<\infty ,$$   see (1.33) \cite{yajima}. Next we look at the
first term in the  rhs of \eqref{Katosquareroot}. We have
\begin{equation}
\label{Laplace2} \begin{aligned} & (\mathcal{H} +m^2+\tau )^{-1}
=(\uno + ( -\Delta  +m^2+\tau )^{-1} V)  ^{-1} ( -\Delta  +m^2+\tau
)^{-1}.
\end{aligned}  \end{equation}
Both factors in the rhs are for $\tau \in [0,a]$  uniformly bounded
as operators from $L^{2,s}$ to itself. In particular, for the first
this can be shown easily to follow  by $|  V(x)|\le C
  \langle x\rangle ^{-5 }$, by Rellich compactness criterion, by
   Fredholm
  theory and by the fact that $\ker (\mathcal{H} +m^2+\tau )=0$
in $L^{2,s}(\R ^3)$ for all
  $\tau \ge 0$ and $s\ge 0$. Hence,
\begin{equation}
   \| \int _0^a
  (\mathcal{H} +m^2+\tau )^{-1}   \frac{d\tau
}{\tau ^{1-\kappa }} \| _{B(L^{2,s}, L^{2,s})} <\infty . \nonumber
\end{equation}
 \qed

 Claim (2) in Lemma \ref{sol.homo1} is a consequence of the following lemma:
\begin{lemma}\label{lem:regularization2}
Assume that $V$ satisfies (H1).
Then, for $\Phi \in \mathcal{S}(\R ^3, \C )$    and for any $\kappa \in \R $ we have $ B^{2\kappa}\Phi \in \mathcal{S}(\R ^3, \C )$.
\end{lemma}
\proof Let us start with $ \kappa >0$. It is elementary that (H1) implies
 $B^{2l}\Phi \in \mathcal{S} (\R ^3, \C )$
 for all $l\in \mathbb{N}$. So it
  is not restrictive to consider
$\kappa <1$.  Then by
  Lemma \ref{lem:regularization1}
  we have $B^{2l+2\kappa }\Phi \in L^{2,s} (\R ^3, \C )$
  for all $l\in \mathbb{N}$ and $s\ge 0$. By (H1)
   this implies also $(-\Delta +m^2)^{l}B^{ 2\kappa } \Phi
   \in L^{2,s} (\R ^3, \C )$ for all $l\in \mathbb{N}$ and $s\ge 0$.
    Hence $B^{ 2\kappa }\Phi \in \mathcal{S} (\R ^3, \C )$. Case $  \kappa =0$
    is elementary by $B^0=P_c$. For  $ \kappa =-2\ell$  with $\ell \in \mathbb{N}$ we can repeat the above proof using the fact that $
  (\mathcal{H} +m^2  )^{-1}   \in {B(L^{2,s}, L^{2,s})}  $ for any
$s\ge 0$. For  more general $ \kappa <0$ for $[|\kappa |]= \ell \in \Z$  for $\ell \le  |\kappa | <\ell +1$  we write $B^{2\kappa}= B^{-2\ell-2} B^{2\kappa + 2\ell +2}$. Then   $\Psi :=B^{2\kappa + 2\ell +2}\Phi\in \mathcal{S}$ because $2\kappa + 2\ell +2>0$ and $B^{2\kappa}\Phi = B^{-2\ell-2}\Psi \in \mathcal{S}$ because $2\ell +2\in \mathbb{Z}$.
  \qed

\bigskip
\noindent  {\bf Proof of Claim (1) Lemma \ref{sol.homo1}.}  We can write
\begin{equation}
    \frac{1}{B -\lambda }\Phi =\frac{1}{B^2 +\lambda ^2}  \Psi  , \quad \Psi :=\lambda \Phi +   B\Phi  . \nonumber
\end{equation}
Since $\Phi \in \mathcal{S} (\R ^3, \C )$ by hypothesis, then
$\Psi \in \mathcal{S} (\R ^3, \C )$ by Lemma \ref{lem:regularization2}.
 By
repeating
 the   argument Lemma \ref{lem:regularization2}  we conclude that  $\frac{1}{B^2 +\lambda ^2}\Psi \in \mathcal{S} (\R ^3, \C )$. Indeed we have $B^{2l}\frac{1}{B^2 +\lambda ^2}\Psi =  \frac{1}{B^2 +\lambda ^2}B^{2l}\Psi   \in L^{2,s}  $
 for all $l\in \mathbb{N}$ and all $s>0$, and this is equivalent to
 $\frac{1}{B^2 +\lambda ^2}\Psi \in \mathcal{S} (\R ^3, \C )$.
 \qed



Dipartimento di Matematica ``Federico Enriques'', Universit\`a degli
Studi di Milano, Via Saldini 50, 20133 Milano, Italy.

{\it E-mail Address}: {\tt dario.bambusi@unimi.it}

\vskip10pt

DISMI University of Modena and Reggio Emilia, Via Amendola 2, Pad.
Morselli, Reggio Emilia 42122, Italy.

{\it E-mail Address}: {\tt cuccagna.scipio@unimore.it}

\end{document}